\title {      Noise Sensitivity of Boolean Functions
\\               And Applications to Percolation                }
\author{
           Itai Benjamini \and Gil Kalai \and Oded Schramm      }
\date{May 15, 2000}
\documentclass[12pt]{article}
\usepackage{amsmath}
\usepackage{amsthm}
\usepackage{amsfonts}
\numberwithin{equation}{section}

\theoremstyle{plain}
\newtheorem{theo}{Theorem}[section]
\newtheorem{thm}[theo]{Theorem}
\newtheorem{lemma}[theo]{Lemma}
\newtheorem{cor}[theo]{Corollary}
\newtheorem{prop}[theo]{Proposition}

\theoremstyle{remark}
\newtheorem{rem}[theo]{Remark}
\newtheorem{rem2}[theo]{Late Remark}

\newtheorem{que}{Question}[section]
\newtheorem{conj}[que]{Conjecture}
\newtheorem{prob}[que]{Problem}

\def\Bbb#1{\mathbb{#1}}
\newcommand{\R}{{\Bbb R}}
\newcommand{\Z}{{\Bbb Z}}

\newcommand{\N}{{\Bbb N}}

\def\Cal#1{{\cal #1}}

\def\P{{\bf P}}
\def\Pdy{\tilde\P}
\def\mu{\P}
\def\E{{\bf E}}

\def\br{\overline}
\def\maj{\Cal M}
\def\Maj{M}
\def\A{\Cal A}
\def\cro{\Cal C}
\def\C{{\cro}}
\def\nn{[n]}
\def\NN{N}  
\def\tilde{\widetilde}
\def\ff{{\widetilde {f}}}
\def\H{{I\!I}}
\def\eps{\epsilon}
\def\neps{\eta}

\def\le{\leq}
\def\ge{\geq}
\def\setminus{-}
\def\hat{\widehat}
\def\var{{\rm {VAR}}}

\def\sss{Z}
\let\CHI=\chi
\def\chi{\raise2.2pt\hbox{$\CHI$}}

\def\proof{\par\noindent{\sc Proof: }}
\def\proofof#1{\medskip\par\noindent{\sc Proof of \ref{#1}: }}
\def\ep{\hfill{$\Box$}}
\let\ep=\qed

\begin{document}

\maketitle

\begin {abstract}
It is shown that a large class of events in a product
probability space are highly sensitive to noise,
in the sense that with high probability, the configuration
with an arbitrary small percent of random errors gives almost no prediction
whether the event occurs.
On the other hand, weighted majority functions are shown to be noise-stable.
Several necessary and sufficient
conditions for noise sensitivity and stability are given.

Consider, for example, bond percolation on an $n+1$ by $n$ grid.
A configuration is a function that
assigns to every edge the value $0$ or $1$.
Let $\omega$ be a random configuration, selected
according to the uniform measure.
A crossing is a path that joins the left and
right sides of the rectangle, and consists entirely of edges
$e$ with $\omega(e)=1$.
By duality, the probability for having a crossing is $1/2$. 
Fix an $\epsilon\in(0,1)$.
For each edge $e$, let $\omega'(e)=\omega(e)$ with probability
$1-\epsilon$, and $\omega'(e)=1-\omega(e)$ with probability $\epsilon$,
independently of the other edges.
Let $p(\tau)$ be the probability for
having a crossing in $\omega$, conditioned on $\omega'=\tau$.
Then for all $n$ sufficiently large,
$\P\big\{\tau : |p(\tau)-1/2|>\epsilon\big\}<\epsilon$.
\end {abstract}

\newpage
\tableofcontents
\newpage

\section {Introduction} \label{s.intro}
\subsection {Noise sensitivity --- three examples}
Consider the Hamming cube $\Omega_n=\{0,1\}^n$
endowed with the uniform probability measure $\P$.
Let $\A\subset \Omega_n$ be some event.
Given a random $x=(x_1,\dots,x_n)\in \Omega_n$,
suppose that $y=(y_1,\dots,y_n)$ is a random
perturbation of $x$; that is, for every $j\in \{1,\dots,n\}$,
$y_j=x_j$ with probability $1-\epsilon$, independently
for distinct $j$'s. 
Here $\epsilon\in(0,1)$ is some small fixed constant.
This random perturbation of $x$ will be denoted $\NN_\eps(x)$.
We may think of $\NN_\eps(x)$ as $x$ with some noise.

Based on the knowledge of $\NN_\eps(x)$, we would like to predict
the event $x\in \A$.
Since the joint distribution $(x,\NN_\eps(x))$ is
the same as that of $(\NN_\eps(x),x)$, an equivalent
problem is to predict $\NN_\eps(x)\in\A$ knowing $x$. 
The event $\A$ is
{\bf noise sensitive} if for all but a small set of
$x$, knowing $x$ does not significantly help
in predicting the event $\NN_\eps(x)\in \A$.
More formally, $\A$ is noise sensitive, if
for some small $\delta>0$,
\begin{equation}
\label{e.delt}
\gamma(\A,\eps,\delta):=
\P\Bigl\{x:\Bigl| \P\big(\NN_\eps(x)\in \A \mid x\big)-
         \P (\A)\Bigr|>\delta\Bigr\}<\delta
\end{equation}

Set
$$
\phi(\A,\eps) = \inf\Bigl\{\delta>0 : \gamma(\A,\eps,\delta)<\delta\Bigr\},
$$
which is the infimum of all $\delta>0$
such that (\ref{e.delt}) holds.  This will be called
the {\bf sensitivity gauge} of $\A$.
A sequence of events $\A_m\subset\Omega_{n_m}$
will be called {\bf asymptotically noise sensitive} if
$$
\lim_{m\to\infty} \phi(\A_m,\eps) = 0, \qquad \forall \epsilon\in(0,1/2).
$$
\begin {rem}
As shown in Section  \ref{s.noise}, $\A_m$ are asymptotically noise
sensitive 
if and only if 
\begin {equation}
\lim _{m \to \infty} \var[\P (N_\eps(x) \in A_m | x)] = 0.
\end {equation}
\end {rem}

A simple example of a sequence of events which are not noise sensitive
is dictatorship.  The first bit dictator is the
event $\Cal D_n=\Bigl\{(x_1,\dots,x_n)\in\Omega_n:x_1=1\Bigr\}$.
To verify that $\{\Cal D_n\}$
is not asymptotically noise sensitive, consider
some event $\A\subset\Omega_n$.
Then for $k>n$ we may obviously consider $\A$
as a subset of $\Omega_k$, by ignoring the extra variables. 
Note that this does not change the value of $\phi(\A,\epsilon)$.
Consequently,
$\phi(\Cal D_n,\epsilon)=\phi(\Cal D_1,\epsilon)\neq 0$
for all $n>1$.
 
Let us examine now the example of majority.
Pick some $\epsilon\in(0,1/2)$.
Let $\maj_n\subset \Omega_n$ denote the majority event,
that is,
$$
\maj_n = \left\{(x_1,\dots,x_n)\in \Omega_n : \sum_j x_j \geq n/2\right\}
.
$$
The probability that $\sum_j x_j-n/2 >\sqrt n$
is bounded from below as $n\to\infty$.
Given such an $x$, the probability that
$\NN_\eps(x)\in\maj_n$ is greater than $\P[\maj_n]+\delta_1$
for some constant $\delta_1>0$, depending on $\epsilon$.
We conclude that majority is not asymptotically
noise sensitive as $n\to\infty$.

Majority and dictatorship are not only noise insensitive,
they are actually ``noise stable'', 
in a sense defined in Subsection~\ref{s.introstab} below.

It turns out that the noise insensitivity of majority
and dictatorship
is atypical, and many natural and interesting events
are asymptotically noise sensitive.

Our third example is bond percolation on an $m+1$ by $m$ rectangle
in the ordinary square grid $\Z^2$.
A configuration is an element in $\Omega=\{0,1\}^E$,
where $E$ is the set of edges in this rectangle.
Let $\omega\in\Omega$ be a random configuration, selected
according to the uniform measure.
A {\bf crossing} is a path that joins the left and
right sides of the rectangle, and consists entirely of edges
$e$ with $\omega(e)=1$.
Let $\cro_m$ be the event that there is some crossing
of this rectangle.
By duality, it is not hard to
see that $\P[\cro_m]=1/2$.

\begin{thm}
\label{t.perc}
The crossing events $\cro_m$ are asymptotically noise sensitive;
that is, $\phi(\cro_m,\epsilon)\to 0$ as $m\to\infty$.
\end{thm}

This theorem will appear as a corollary of a 
general result.  To introduce the more general
statement, we need the notion of influence.

\subsection {Influences of variables}
\medskip

Set $[n]=\{1,\dots,n\}$.
Given $x\in \Omega$ and $j\in\nn$, let
$\sigma_jx=(x'_1,\dots,x'_n)$, where
$x'_k=x_k$ when $k\neq j$ and $x'_j=1-x_j$.
The {\bf influence}  of  the $k$-th 
variable on a function $f:\Omega\to\R$ is defined by
\begin{equation}\label{e.idef}
I_k(f) = \| f(\sigma_kx)-f(x)\|_1
\,.
\end{equation}
In other words, $I_k(f)$ is the expected absolute value
of the change in $f$ when the $k$'th bit $x_k$ is flipped.
We shall often not distinguish between an event
$\A$ and its indicator function $\chi_{\A}$.
In particular, for events $\A$, $I_k(\A)=I_k(\chi_\A)$.
Note that $I_k({\A})$ is the probability that
precisely one of the two elements $x,\sigma_kx$ is in ${\A}$.

This notion of influence was introduced by Ben-Or and Linial~\cite {BL}.
Kahn, Kalai and Linial \cite {KKL} (see also, \cite {BKKKL,T1})
showed that
for every $\A\subset\Omega_n$ with $\P[\A]=1/2$
there is a $j\in\nn$ with $I_j(\A)\geq c \log n /n$, for some
constant $c>0$,
and that there always exists
a set $S\subset\nn$ with $|S|\leq c(\epsilon)n/ \log n$
whose cumulative influence is
$> 1-\epsilon$; that is, the measure of the set of inputs for variables
in $[n]\setminus S$ which determine the value of $f$ is less than
$\epsilon$.  

Put 
\begin{eqnarray*}
I(f)&=&\sum_k I_k(f),
\\
\H(f) &=& \sum_k I_k(f)^2
.
\end{eqnarray*}

\begin{thm}
\label{main}
Let $\A_m\subset\Omega_{n_m}$ be a sequence of events
and suppose that $\H(\A_m)\to 0$ as $m\to\infty$.
Then $\{\A_m\}$ is asymptotically noise sensitive.

Equivalently, there is some continuous function
$\Phi$ satisfying $\Phi(0,\epsilon)=0$ such that
$\phi(\A,\eps)\leq \Phi\big(\H(\A),\eps\big)$ for
every event $\A$ in some $\Omega_n$.
\end{thm}

On $\Omega_n$, we use the usual lattice order:
$(x_1,\dots,x_n)\leq(y_1,\dots,y_n)$ iff
$x_j\leq y_j$ for all $j\in\nn$.
A function $f:\Omega_n\to\R$ is {\bf monotone}
if $f(x)\leq f(y)$ whenever $x\leq y$.
An event $\A\subset\Omega_n$ is {\bf monotone} if
its indicator function $\chi_{\A}$ is monotone.

For monotone events, Theorem~\ref{main} has a converse:

\begin{thm}
\label{p.mconverse}
Let $\A_m\subset\Omega_{n_m}$ be a sequence of monotone events
with 
$$
\inf_m \H(\A_m) >0
.
$$
Then $\{\A_m\}$ is not asymptotically noise sensitive.
\end{thm}

The assumption that the events $\A_m$ are monotone is
necessary here.  (For example, take $\A_m$ to be
a uniform random subset of $\Omega_{m}$, or parity:
$\A_m:=\{x\in\Omega_m :\, \|x\|_1\hbox{ is odd}\}$.)

\medskip 
Suppose that $\A$ is a monotone event where the influences of all the
variables
are the same.  The influence $I_1(\A)$ then measures the sensitivity of 
$\A$ to flips of a single variable.
Note that, quite paradoxically, $\A$ is {\it least\/}
sensitive to noise when $I_1(\A)$ is largest.

\medskip 
We now give a quantitative version of Theorem~\ref{main} under the
assumption that $\H(\A_m)$ goes to zero fast enough.

\begin{thm}
\label{t.pow}
Let $\A\subset\Omega_{n}$, and suppose that
$\H(\A)\leq n^{-a}$, where $a\in(0,1/2]$.
Then there exist  $c_1,c_2>0$, depending only on $a$ so that
$$
\phi(\A,\eps)\leq c_1 n^{-c_2 \eps}
\,,\qquad \forall \eps\in(0,1/4)
\,.
$$
Consequently, if $\A_m\subset\Omega_{n_m}$
is a sequence  of events satisfying $\H(\A_m)\leq (n_m)^{-a}$ and $\eps_m$
is a sequence in
$(0,1/4)$ such that $\eps_m\log n_m\to\infty$, then $\phi(\A_m,\eps_m)\to
0$.
\end{thm}

\subsection {Weighted majority}
It turns out that for monotone events noise insensitivity is
also closely related to correlation with majority functions.

Let $K\subset\nn$ and define the {\bf majority function} on $K$ by
$\Maj_K(x) = \hbox{\rm sign}\sum_{j\in K} (2 x_j-1)$; that is,
\begin{equation}\label{e.majdef}
\Maj_K(x)=
\left\{
\begin{array}{rl}
{}-1 & \mbox{if $\sum_{j\in K} x_j < |K|/2$}\,;\\
0 & \mbox{if $\sum_{j\in K} x_j = |K|/2$}\,;\\
1 & \mbox{if $\sum_{j\in K} x_j > |K|/2$}\,.
\end{array}\right.
\end{equation}
For $f:\Omega_n\to\R$ set
$$
\Lambda(f)=\max\Bigl\{|\E(f\Maj_K)|:K\subset\nn\Bigr\}
.
$$

\begin{thm}\label{t.lambii}
Let $f:\Omega_n\to [0,1]$ be monotone.
Then
$$
\H(f)\leq C \Lambda(f)^2\Bigl(1-\log\Lambda(f)\Bigr)\log n
,
$$
where $C$ is some universal constant.

Consequently, if
$\A_m\subset\Omega_{n_m}$ is a sequence of monotone events
with
\begin{equation}\label{e.majcrit}
\lim_{m\to\infty}
\Lambda(\A_m)^2\Bigl(1-\log\Lambda(\A_m)\Bigr)\log n_m = 0
.
\end{equation}
Then $\{\A_m\}$ is asymptotically noise sensitive.
\end{thm}

One cannot get rid of the $\log n_m$ factor (see Remark~\ref{r.nouni}),
except by using weighted majority functions.
For positive weights ${\bf w}= (w_1,w_2, \dots, w_n)$ 
consider a weighted majority function, which is 
defined  by
$$
M_{\bf w}(x_1, x_2, \dots, x_n)= \hbox{\rm sign}\left(\sum (2 x_j-1)w_j
\right)
.
$$
Finally write 
$$\tilde\Lambda (A)=
\max\Bigl\{|\E(\chi_A M_{\bf w})|:
{\bf w} \in[0,1]^n\Bigr\}.$$

\begin{thm}\label{t.lambii2}
Let $\A_m\subset\Omega_{n_m}$ be a sequence of monotone events.
Then $\{\A_m\}$ is asymptotically noise sensitive if and only if 
$\lim_{m\to\infty} \tilde\Lambda (A_m) =0$. 
\end{thm}

For a monotone event
$\A\subset\Omega_n$, which is symmetric in the $n$ variables,
its correlation with unweighted majority is enough to determine
if it is noise sensitive.

\subsection {Stability}\label{s.introstab}

We now define the notion of {\bf stability}, which is the opposite of
noise sensitivity.  Suppose $\A\subset\Omega_n$, and let
$x\in\Omega_n$ be random-uniform.
For $\eps>0$, let $\NN_\eps\A$ denote the event $\NN_\eps(x)\in\A$.
It is then clear that $\P[\A\triangle\NN_\eps\A]\to 0$ as $\eps\to0$.
($\Cal B\triangle\A$ denotes the symmetric difference,
$(\Cal B-\A)\cup(\A-\Cal B)$.)  The faster 
$\P[\A\triangle\NN_\eps\A]$ tends to zero, the more noise-stable
$\A$ is.  
More precisely, let $\{\A_i \}$ be a collection
of events, where $\A_i\subset\Omega_{n_i}$.  We say that
$\{\A_i\}$ are {\bf uniformly stable} if the limit
$\lim_{\eps\to0}\P[x \in \A_i\triangle\NN_\eps\A_i]=0$ is uniform in $i$.

For $\bf w\in\R^n$ and $s\in\R$, let $\Cal M_{{\bf w},s}$ be the
(generalized)
weighted majority event
$$
\Cal M_{{\bf w},s}:=
\left\{x\in\R^n: \sum_{j=1}^n (2x_j-1)w_j>s\right\}\subset\Omega_n.
$$
Let $\mathfrak M$ denote the collection of such events:
$$
\mathfrak M :=\bigl\{\Cal M_{{\bf w},s}:n=1,2,\dots,\ \bf w\in\R^n,\
s\in\R\Bigr\}
.
$$
In Section \ref{s.maj} we show that 

\begin{thm}\label{t.stab}
$\mathfrak M$ is uniformly stable.  Moreover, for every $\Cal M\in\mathfrak
M$
$$
\P[\Cal M - \NN_\eps\Cal M]\leq C \eps^{1/4},
$$
where $C$ is a universal constant independent of $\Cal M$.
\end{thm}

Note that an infinite sequence $\{\A_i\}$ with
$\P[\A_i]$ bounded away from $0$ and $1$ cannot be asymptotically
noise sensitive and uniformly stable. 
We also observe (Lemma~\ref{l.corstab}) that
when $\{\A_i\}$, ($\A_i\subset\Omega_{n_i}$),
is asymptotically noise sensitive and
$\{\Cal B_i\}$, ($\Cal B_i\subset\Omega_{n_i}$),
is uniformly stable, then $\A_i$ and $\Cal B_i$ are asymptotically
uncorrelated.  One can say, somewhat imprecisely, that
the noise sensitive events are asymptotically in the
orthocomplement of the uniformly stable events.

Stability and sensitivity are two extremes.  However,
there are events that are neither sensitive nor stable.
For example, if $\Cal C$ is the event of a percolation
crossing, as described above, and $\Cal M$ is the majority
event, then $\Cal C\cap \Cal M$ is neither asymptotically
noise sensitive, nor uniformly stable.

\subsection {Fourier-Walsh expansion}

For a boolean function $f$ on $\{0,1\}^n$,
consider the Fourier-Walsh expansion 
$ f = \sum_{S \subset [n]} \hat{f}(S) u_S,$ 
where, $u_S(T)=(-1)^{|S \cap T|}$.
Here and in the following, we identify any vector 
$x\in\Omega^n$ with the subset $\{j\in[n]:x_j=1\}$, of
$[n]=\{1,2,\dots,n\}$.
Consequently, $|x|$ denotes the cardinality of that set; that is,
$|x|=\|x\|_1$
for $x\in\Omega_n$.

\begin{thm}\label{t.foursens}
Let $\A_m\subset\Omega_{n_m}$ be a sequence of events, and
set $g_m=\chi_{\A_m}$.
Then $\{\A_m\}$ is asymptotically noise sensitive iff for every
finite $k$
\begin{equation}\label{e.foursens}
\lim_m
\sum\Bigl\{\hat g_m(S)^2 : S\subset\nn, 1\leq |S|\leq k
\Bigr\} = 0
.
\end{equation}
$\{\A_m\}$ is uniformly stable iff 
\begin{equation}\label{e.fourstab}
\lim_{k\to\infty}\sup_m
\sum\Bigl\{\hat g_m(S)^2 : S\subset\nn, |S|\ge k
\Bigr\} = 0
.
\end{equation}
\end{thm}

It can be easily shown that  for $f=\chi_{\A}$
$$I (f) = 4 \sum_{S \subset [n]} \hat f(S)^2|S|. $$
(This follows from~(\ref{e.fjnrmp}) below with $p=2$.)
We will introduce another quantity
$$J (f) = \sum_{\emptyset\neq S \subset [n]} \hat f(S)^2/|S|. $$
Also set for $\A\subset\Omega_n$, $n>1$,
\begin{align*}
\alpha (\A) &= \log I(\A)/\log n,
\\
\beta (\A) &= -\log J(\A) /\log n.
\end{align*}
For events $\A$ we clearly have $0\le\beta(\A)$,
and $\beta(\A)\le\alpha(\A)$, provided that $\P[\A]= 1/2$.
When $\A$ is monotone $\alpha (\A) \le 1/2$.

Perhaps some words of explanation are needed. $I(\A)$ measures the 
sum of the influences of the variables. 
For monotone events
it is maximal for majority, where $I(\A)\simeq \sqrt n$
and thus $\alpha (\A)\to 1/2$. 
In the terminology
used in percolation theory, $I(\A)$ is the expected number of pivotal edges.

For the crossing events $\Cal C$ of percolation (in arbitrary dimensions) 
it is conjectured that $I(\Cal C)$ behaves like a certain fractional power 
(a {\it critical exponent\/}) of $n$.
It is conjectured that in dimension 2,
as $n$ tends to infinity, $\alpha (\Cal C)$ tends to 3/8.
Thus, this critical exponent generalizes and has a 
Fourier-analysis interpretation for arbitrary Boolean functions.

$\alpha(\A)$ is large if there are substantial Fourier coefficients
$\hat f(S)$ for large $|S|$. In contrast,  
$\beta(\A)$ is large if there are {\it no\/} substantial 
Fourier coefficients $\hat f(S)$ for $S$ of small positive size.
We conjecture that for the crossing events for 
percolation, as $n$ tends to infinity  $\beta (\C)$ 
tends to a positive limit. We are curious to know whether this 
limit is strictly smaller than 
the limit for $\alpha(\C)$.

\subsection {Some related and future work}

There are interesting connections between noise sensitivity and 
isoperimetric inequalities of the form described by
Talagrand in  \cite {T1.5}. These connections 
and applications for first passage
percolation problems will be discussed in a subsequent paper, \cite {BKS2}.

Our notion of noise sensitivity is related to the study 
of noises by Tsirelson~\cite {Ts1, Ts2}. ``Noise'', in Tsirelson's sense, is
a
type of $\sigma$-field filtration. Uniform stability 
seems to correspond, in the limit,
to the noise being white, while asymptotic sensitivity 
seems to correspond to the noise being black.

\subsection {The structure of this paper}

Theorems \ref{main}, and \ref{p.mconverse}
are proved in the next section. 
Our proofs combines combinatorial reasonings with applying 
certain inequalities for the Fourier coefficients 
of Bonami and Beckner which were used
already in \cite {KKL}. However, to get the results
in the sharpest forms we have to rely on a sophisticated
\lq\lq bootstrap" method of \cite {T2} and on the main results
of that paper which rely on this method. 
Talagrand's remarkable paper \cite {T2}
has thus much influence on the present work.

Weighted and unweighted majority 
functions are considered in Section~\ref{s.maj}. 
An applications to percolation
is described in Section~\ref{s.perc} followed by some related open problems
in Section~\ref{s.percconj}.
In Section~\ref{s.ex}, we will work out two examples (due
to Ben-Or and Linial). In one of these $\alpha (\A) \to 1-\log_23$ and
$\beta (\A)\to 1- \log_23$.
In Section~\ref{s.compl} we consider relations with complexity theory. 
A simple description of  
noise-sensitivity in terms of random walks 
is given in Section~\ref{s.rw}.
In Section~\ref{s.fix} we consider perturbations with a different
sort of noise, where the number of bits that are changed
is fixed. The conclusions are similar to those above, but
there is an amusing and slightly unexpected twist.

For simplicity we consider here the uniform measure on $\Omega_n$.
More generally, one may consider the product measure $\P_p$,
where $\P_p\{x:x_j=1\}=p$.  Our results and proof apply
in this setting.
(All that is needed is to replace the Fourier-Walsh
transform by its analog as given in Talagrand's paper \cite {T1}
and the proofs go through without change.)
However, the case when $p$ itself depends on $n$
is interesting, but will not be considered here.

Since the first version of this paper was distributed, a few of the 
problems we posed were settled by several people, not always in the 
direction anticipated by us. These developments are mentioned briefly 
in a few ``late remarks'' throughout the paper.

\medskip

{\bf\noindent Acknowledgments.}

It is a pleasure to thank Noga Alon, 
Ehud Friedgut, Ravi Kannan, 
Harry Kesten, Yuval Peres, Michel Talagrand and 
Avi Wigderson for helpful discussions. 

\section {Sensitivity to noise} \label{s.noise}

We now put the noise operator $\NN_\eps$ defined
in the introduction into a somewhat more general framework.
That will allow us to deal, for example, with the situation
where the $1$ bits are immune to noise but the $0$ bits
are noise prone.

Consider the following method for 
selecting a random point $x\in \Omega_n$.
Let $q_1,\dots,q_n$ be independent random variables in $[0,1]$,
with $\E q_j=1/2$, for $j=1,\dots,n$, and let
$\omega\in[0,1]^n$ be random uniform.
Set 
$$
x_j=
\left\{\begin{array}{ll}
1, & \mbox{if }1-\omega_j<q_j,\\
0,&\mbox{otherwise}.
\end{array}
\right.
$$
Then $x$ is distributed according to the uniform measure of $\Omega_n$;
it will be denoted by $\NN(\omega,q)$.

Let $\nu$ be the measure on $[0,1]^n$ such that
$\nu(X)=\P\bigl[(q_1,\dots,q_n)\in X\bigr]$.
We think of $x$ as being chosen in two stages.
In the first stage, $q=(q_1,\dots,q_n)$
is selected according to $\nu$.  This $q$ gives a product measure $\P_q$
on $\{0,1\}^n$ that satisfies
$\P_q\{\tau\in \Omega_n : \tau(j)=1\} = q_j$.
Then $x$ is chosen according to the measure $\P_q$.

For example, suppose $z\in\Omega_n$.  Define $q=q(z)\in[0,1]^n$ by
$q_j=1-\eps$ if $z_j=1$ and $q_j=\eps$ if $z_j=0$.
Then for every $z\in\Omega_n$, the perturbation $\NN_\eps(z)$
has the same distribution as $\NN(\omega,q(z))$.
The $\nu$ giving this distribution of $q$ will be denoted
$\nu_\eps$.

However, the construction $\NN(\omega,q)$
is more general than that given by the noise operator $\NN_\eps$.
As hinted above, one can create a situation where $1$ bits are
robust, but $0$ bits are prone to noise.  More precisely,
take $q_j=1$, with probability $1/2-\epsilon$ and
$q_j=\epsilon/(1/2+\epsilon)$ with probability $1/2+\epsilon$.

Another interesting example is obtained by taking each $q_j$
to be $1$, with probability $(1-\epsilon)/2$,
$0$, with probability $(1-\epsilon)/2$,
and $1/2$ with probability $\epsilon$.

Let $f:\{0,1\}^n\to \R$ be some function.  In the following,
$f$ will be taken to be the characteristic function $\chi_\A$ of
some event $\A\subset\{0,1\}^n$, or $f=\chi_\A-\mu(\A)$.
What information does the first stage in the selection
of $x=\NN(\omega,q)$, namely the selection of $q$, give about
the value of $f(x)$?
If we know that $q=z$, then our prediction for $f(x)$
would be 
$$
G(f,z)=\E\big(f(x) \mid q=z\big).
$$
The expected value of $G(f,q)$ is obviously $\E(f)$.
Let
$$
\sss(f,\nu) = \E_q G(f,q)^2 = \int G(f,z)^2 \,d\nu(z) 
. 
$$
This is just the second moment of $G(f,q)$. 
If $\sss(f,\nu)-(\E f)^2$
is small, then for ``most'' values of
$q$ there is no prediction for $f(x)$ that
is significantly better than the a~priori knowledge of $\E f$.
We often write $G(\A,\cdot)$ and $\sss(\A,\cdot)$
in place of $G(\chi_\A,\cdot)$ and $\sss(\chi_\A,\cdot)$.

\begin {lemma}\label{l.gauge}
The number $\sss(f,\nu)$ depends only
on $f$ and the variances $\zeta_j$ of the variables $q_j$.
Its expression in terms of the Fourier coefficients is,
$$
\sss(f,\nu) = \sum_{S\in\Omega_n}\hat f(S)^2 \prod_{j\in S} 4 \zeta_j
.
$$
\end {lemma}

\proof
\begin{align*}
G(f,z) &= \E\big(f(x)\mid q=z\big)
\\ &= 
\sum_{T\subset[n]} f(T) \prod_{j\in T} z_j\prod_{j\notin T}(1-z_j)
\\ &=
\sum_ T\sum_S
\hat f(S) (-1)^{| T\cap S|}  \prod_{j\in T} z_j\prod_{j\notin T}(1-z_j)
\\ &=
\sum_S \hat f(S)
\left( \sum_{ T'\subset S} (-1)^{| T'|} 
\prod_{j\in T'} z_j\prod_{j\notin T'}(1-z_j)\right)
\cdot \\&\qquad\qquad\qquad\qquad
\left( \sum_{ T''\subset[n]-S}\,\, 
\prod_{j\in T''} z_j\prod_{j\notin T''}(1-z_j)\right)
\\ &=
\sum_S \hat f(S)
\left(\prod_{j\in S} \big((1-z_j)-z_j)\right)
\left(\prod_{j\notin S} \big((1-z_j)+z_j)\right)
\\ &=
\sum_S \hat f(S)
\prod_{j\in S} (1-2z_j)
\end{align*}
Therefore,
\begin{eqnarray*}
\sss(f,\nu)
&=&
\E G(f,q)^2
\\&=&
\sum_S\sum_{S'} \hat f(S)\hat f(S')
\E\left(\prod_{j\in S} (1-2 q_j)\prod_{j\in S'}(1-2 q_j)\right)
\\&=&
\sum_S\sum_{S'} \hat f(S)\hat f(S')
\prod_{j\in S\cap S'} \E(1-2 q_j)^2
\prod_{j\in S\triangle S'} \E(1-2 q_j)
.
\end{eqnarray*}
Since $\E q_j=1/2$, summands with $S\neq S'$ vanish.
The lemma follows.
\ep \medskip

For every $\eps\in[0,1]$, $x\in\Omega_m$ and $f:\Omega_n\to\R$ set
$$
Q_\eps f(x) = \E f(\NN_\eps(x))
$$
(here the expectation is only with respect to the noise).
Also let
$$
\var(f,\epsilon) =
{\rm var}(Q_\eps f) = \sss\big(f,\nu_\eps\big) - (\E f)^2
\,.
$$
Note that for singletons $S=\{i\}\subset[n]$, we have
$Q_\eps u_S = (1-2\eps) u_S$.
If $S_1,S_2\subset[n]$ are disjoint and $x\in\Omega_n$ is fixed,
then $N_\eps(x)\cap S_1$ and $N_\eps(x)\cap S_2$ are independent.
Consequently, $Q_\eps (u_{S_1}u_{S_2}) = (Q_\eps u_{S_1})(Q_\eps u_{S_2})$.
We may conclude that $Q_\eps u_S = (1-2\eps)^{|S|} u_S$ for
every $S\subset[n]$, and linearity gives 
\begin{equation}
\label{e.qq} 
Q_\eps f = \sum_{S\subset[n]} \hat f(S)(1-2\eps)^{|S|} u_S
\,.
\end{equation}
One consequence of this, which can also be obtained from
Lemma~\ref{l.gauge}, is
\begin{equation}\label{e.vareps}
\var(f,\epsilon) = 
\sum_{\emptyset\neq S\subset [n]}
\hat f(S)^2 (1-2\epsilon)^{2|S|}
.
\end{equation}

Now we relate $\var(\A,\eps)$ with the sensitivity gauge
$\phi(\A,\eps)$:

\begin{prop}\label{p.varphi}
For every $\A\subset\Omega_n$
$$
\frac{1}{2}\var(\A,\eps)\leq \phi(\A,\eps) \leq \var(\A,\eps)^{1/3}
.
$$
\end{prop}

\proof
Let $\delta=\phi(\A,\eps)$, and set
$$
Y=\Biggl\{y\in\Omega_n:\Bigl|Q_\eps\chi_\A(y) - \P[\A]\Bigr|\geq
\delta\Biggr\}
.
$$
Then, by the definition of $\phi$, $\P[Y]\geq\delta$.
Consequently,
$$
\var(\A,\eps)
\geq\delta^2\P[Y]\geq\delta^3=\phi(\A,\eps)^3
.
$$

For the other direction set
$$
Y'=\Biggl\{y\in\Omega_n:\Bigl|Q_\eps\chi_\A(y) - \P[\A]\Bigr|>
\delta\Biggr\}
.
$$
Then $\P[Y']\leq \delta$.
For $y\in Y'$, the trivial estimate
$\Bigl|Q_\eps\chi_\A - \P[\A]\Bigr|\leq 1$ holds.
Therefore,
$$
\var(\A,\eps)
\leq \P[Y'] + \delta^2\leq 2 \phi(\A,\eps)
.
$$
\ep
\medskip

\proofof{t.foursens}
The first part is immediate from Prop.~\ref{p.varphi} and~(\ref{e.vareps}).
For the proof of the second part, observe that~(\ref{e.qq}) implies
that~(\ref{e.fourstab}) is equivalent to $\|g_m-Q_\eps g_m\|_2\to 0$
uniformly
as $\eps\to 0$.
Since $|g_m|$ and $|Q_\eps g_m|$ are bounded, this is equivalent
to $\|g_m-Q_\eps g_m\|_1\to 0$ uniformly, which is the same as
uniform stability for $\{\A_m\}$.
\ep
\medskip 

\begin{rem}
\label{r.enough}
Another consequence of~\ref{p.varphi} and~(\ref{e.vareps}) is that
for constant $\eps,\eps'\in(0,1/2)$, we have
$\phi(\A_m,\eps)\to 0$ iff $\phi(\A_m,\eps')\to 0$.
Consequently, to verify that $\A_m$ is 
asymptotically noise sensitive, it is enough to prove 
$\var(\A_m,\eps)\to 0$ with any fixed $\eps\in(0,1/2)$.
\end{rem}

By Theorem~\ref{t.foursens}, to establish Theorem~\ref{main} we need to show
that the $L^2$ weight of the Fourier coefficients with
$|S|$  small is negligible.
For a function $g= \sum \hat g (S) u_S$ let 
$$
T_{\neps}g = Q_{\frac{1-\neps}{2}} g =  \sum \hat g(S) \neps ^{|S|}u_S
.
$$
Observe that $T_0(g)=\E g$ and $T_1g=g$.
Also note that
\begin{equation}\label{e.tvar}
\|T_{1-2\eps}g\|^2_2 = \var(g,\eps)+\hat g(\emptyset)^2
,
\end{equation}
by (\ref{e.vareps}).

The following hyper-contractive inequality of Bonami and Beckner~\cite
{Bo,Be},
which was crucial in \cite {KKL}, will be useful.

\begin {lemma} [Bonami, Beckner]
\label{l.bb}
$\|T_{\neps}f\|_2 \le \|f\|_{1+\neps^2}.$
\end {lemma}

The following is a slightly weaker version of Theorem~\ref{main},
which is sufficient for the applications to percolation.
It is presented here, since we can give an almost self-contained proof of
it.

\begin {thm}\label{t.weak}
Suppose that $\A_m\subset\Omega_{n_m}$ is a sequence
of 
events and 
\begin{equation}
\lim_{m\to\infty} {\log\H(\A_m)\over \log\log n_m} = -\infty
\,.
\label{e.fast}
\end{equation}
Then $\{\A_m\}$ is asymptotically noise sensitive.
\end {thm}

\proof
Abbreviate $\A$ for $\A_m$ and $n$ for $n_m$, and set
$f:=\chi_\A$.
Let $\ff:=\chi _\A - \P [\A]$.
Thus, $\hat \ff (\emptyset)=0$ and $\hat \ff (S)=\hat f(S)$,
when $S\neq\emptyset$. 

Recall that $\sigma_jx=(x'_1,\dots,x'_n)$,
where $x'_i=x_i$ if $i\neq j$ and $x'_j=1-x_j$. 
Let 
$$
f_j(x) = f(x) - f(\sigma_j x)
,\qquad j=1,2,\dots,n,
$$
and note that
$$
\hat f_j(S) =
\left\{\begin{array}{ll}
0, & \mbox{if }j\notin S,\\
2\hat f(S),&\mbox{if }j\in S.
\end{array}
\right.
$$
Since $f_j$ takes only the values $-1,0,1$,
equation~(\ref{e.idef}) gives for every $p\geq 1$,
\begin{equation}\label{e.fjnrmp}
\|f_j\|_p =  I_j(f)^{1/p}
\,.
\end{equation}
We set $\neps :=1-2\eps$, where $\eps\in(0,1/2)$ and
$$
F_\A(\neps) := \var(\A,\eps)=
\|T_{\neps}\ff \|_2^2 =
\sum _{S\neq\emptyset}\hat f(S)^2 \neps ^{2|S|} 
.
$$
By Remark~\ref{r.enough} and Theorem~\ref{p.varphi},
it is enough to prove that
$F_\A(1/2)\to 0$ as $m\to\infty$.  We have
\begin{eqnarray}
F_\A(\neps) &\le &
\sum _{S}\hat f(S)^2 |S| \neps ^{2|S|} = 
{1\over4}\sum _{j=1}^n \|T_\neps f_j\|_2^2 
\nonumber
\\ &\le &
{1\over4}\sum _{j=1}^n \|f_j\|_{1+\neps^2}^2 
\qquad\mbox{(by Lemma~\ref{l.bb})}
\nonumber
\\ &\le &
\sum _{j=1}^n I_j(\A)^{2/(1+\neps^2)}
\qquad\mbox{(by (\ref{e.fjnrmp}))}
\nonumber
\\ &\le &
n^{\neps^2/(1+\neps^2)} \H(\A)^{1/(1+\neps^2)}
\qquad\mbox{(by the means inequality)}
.
\label{e.nh}
\end{eqnarray}
Take some $\neps_1\in(0,1/2)$, to be later specified, and set
$\lambda:=\log F_\A(\neps_1)/\log\neps_1$.  If $\neps\geq\neps_1$, then
\begin{eqnarray}
F_\A(\neps)
 & \le &
\sum_{1\leq |S|\leq\lambda/2  } \hat f(S)^2 \neps^{2 |S|}
+ \neps^\lambda \sum_S \hat f(S)^2
\nonumber
\\ & \le &
(\neps/\neps_1)^\lambda F_\A(\neps_1)+\neps^\lambda
= 2 \neps^\lambda
\,.
\label{e.nl}
\end{eqnarray}
Assume that $\H(\A)\in\left(0,e^{-2}\right)$, and let
$a:=\min\Bigl\{-\log\H(\A)/\log n,1/2\Bigr\}$.
We may choose $\neps_1:= \sqrt{a/2}$.
Then $\H(\A)\leq n^{-a}$, and therefore~(\ref{e.nh}) and the
definition of $\lambda$ give
\begin{equation}
\lambda \ge { a\log n\over 3\log(1/a)}
\label{e.lb}
\,.
\end{equation}
The definition of $a$ together with~(\ref{e.fast})
and~(\ref{e.lb}) show that $\lambda\to\infty$ as $m\to\infty$.
Hence~(\ref{e.nl}) implies $F_{\A_m}(1/2)\to 0$ as $m\to\infty$,
which completes the proof.
\ep
\medskip

\proofof{t.pow}  The above calculations together with
Prop.~\ref{p.varphi} show that
$$
\phi(\A,\eps) \leq
\var(\A,\eps)^{1/3} =
F_\A(1-2\eps)^{1/3} \leq 2^{1/3} \Bigl(1-2\eps\Bigr) ^{a\log n\over
9\log(1/a)}
\,,
$$
for $\eps\in(0,1/4)$, when we assume $\H(\A)\leq n^{-a}$, $a\in(0,1/2]$.
The theorem follows immediately.
\ep
\medskip

For the proof of Theorem~\ref{main}, we will need the following.

\begin {thm}
\label{t.tala}
For each $k=1,2,\dots$,
there is a constant $C_k<\infty$ with the following property.
Let $\A\subset\Omega_n$ be a monotone event and $f=\chi_\A$.
Then 
$$
\sum_{|S|=k}  \hat f (S)^2 \le C_k \H(A) \Bigl(-\log \H(A)\Bigr)^{k-1}
\,.
$$ 
\end {thm}

This inequality was proved by Talagrand~\cite {T2} for 
$k=2$. (Talagrand considers an extension of this relation for two events,
and our generalization applies for that extension as well.) 

\medskip
\proofof{t.tala}
To prove the theorem one can follow Talagrand's proof almost word-by-word. 
We will only describe the changes needed to adapt the proof.
One modification required is that the inequality
\begin{equation}
\label{e.subgsub}
\mu \left\{S' :\left| \sum_{|S|=k} \alpha_S u_S(S')\right| \ge t\right\} 
\le e^2 \exp \Biggl(- t ^{2/k} \Bigl(e^2 \sum \alpha
_S^2\Bigr)^{-1/k}\Biggr)\,
\end{equation}
must be used in place of the sub-Gaussian estimate that appears
as Prop.~2.1 in~\cite{T2}.
Set $q=t^{2/k}/\left(e^2 \sum \alpha_S^2\right)^{1/k}$.
For $q\leq 2$ the inequality~(\ref{e.subgsub}) is trivial,
while for $q>2$ it follows by substituting $q$ into
\begin{equation}\label{e.dualBon}
\left\| \sum_{|S|=k} \alpha _S u_S\right\|_q 
\le (q-1)^{k/2}\left(\sum \alpha_S^2\right)^{1/2}\,,
\qquad\forall q\geq 2\,,
\end{equation}
which appears in~\cite {T1} as (2.4) and is a 
consequence of the dual version of the Bonami-Beckner inequality.

Set $A_k:=\{x\in\A: \sigma_k x\notin\A\}$, and note that
$2\P[A_k]=I_k(\A)$.
In the proof for the case $k=2$, Talagrand considers
in Section 3 of~\cite{T2} partitions
$I\cup J=\nn$, and estimates $\sum\{ I_j(\A)^2:{j\in L(s)}\}$,
where $L(s)$ is the set of $j\in J$ such that
$$
\sum_{i\in I}\left(\int_{A_j}u_{\{i\}}(x)\right)^2 \geq s^2 \P[A_j]^2
\,.
$$
To generalize Talagrand's argument for $k>2$,
one gives a similar estimate to $\sum\{ I_j(\A)^2 : {j\in L_{k-1}(s)}\}$,
where $L_{k-1}(s)$ is the set of $j\in J$ such that
$$
\sum\left\{\left(\int_{A_j}u_i(x)\right)^2 :
  i\subset I,\, |i|=k-1\right\}\geq s^2 \P[A_j]^2
$$
We omit the details, since from this point on
only straightforward changes are required to adapt
Talagrand's beautiful (but rather mysterious)
argument.
\ep
\medskip

In the case of monotone events, Theorem~\ref{main}
follows immediately from Theorems~\ref{t.tala}
and~\ref{t.foursens}.
In order to get rid of the monotonicity assumption,
we introduce the shifting operator.

Let $j\in\{1,\dots,n\}$, and let $f:\Omega_n\to\R$.
For $x\in\Omega_n$, set
$$
\kappa_j f(x):=
\left\{\begin{array}{rl}
\max\{f(x),f(\sigma_j x)\}, & \mbox{if }x_j=1,\\
\min\{f(x),f(\sigma_j x)\}, & \mbox{if }x_j=0.\\
\end{array}
\right.
$$
The operator $\kappa_j$ is called the {\bf $j$-shift}.
The following lemma describes some useful properties
of shifts.

\begin{lemma}[Shifting]
\label{l.shift}
Let $f:\Omega_n\to \R$, and let $j,i\in\{1,\dots,n\}$.
Then
\begin{enumerate}
\item $\kappa_1\kappa_2\dots\kappa_n f$ is monotone.
\item $I_{i}(\kappa_j f)\leq I_{i}(f)$.
\item $ \var(\kappa_j f,\eps) \geq \var(f,\eps)$
for each $\eps\in[0,1]$.
\end{enumerate}
\end{lemma}

\proof
Suppose for the moment that $i\neq j$.
For any $a,b\in\{0,1\}$ and
$x\in\Omega_n$, let $x_{a,b}$ be $x$ with the $i$'th coordinate set to
$a$ and the $j$'th coordinate set to $b$.
Note that $\kappa_j f$ is monotone nondecreasing in the variable $x_j$.  
Hence $\kappa_i\kappa_j f(x_{1,1})$ is the maximum
of $f$ on $\{x_{0,0},x_{0,1},x_{1,0},x_{1,1}\}$
and $\kappa_i\kappa_j f(x_{0,0})$ is the minimum.
It follows that $\kappa_i\kappa_j f = \kappa_j\kappa_i\kappa_j f$.
This relation easily implies the first claim of the lemma.

For the second part, we may assume with no loss of generality that $j\neq
i$,
because $I_{i}(\kappa_i f)=I_{i}(f)$. 
A case by case analysis shows that
\begin{eqnarray*}
&&|f(x_{0,0})-f(x_{1,0})| + |f(x_{0,1})-f(x_{1,1})| 
\\
&&\qquad \qquad\geq 
|\kappa_jf(x_{0,0})-\kappa_jf(x_{1,0})| +
|\kappa_jf(x_{0,1})-\kappa_jf(x_{1,1})| 
\,,
\end{eqnarray*}
and the second part follows by summing over $x\in\Omega_n$.

For the last part, set 
\begin{eqnarray*}
g(y) &=& \E[f(\NN_\eps(x))\mid x=y],\\
\tilde g(y) &=& \E[\kappa_j f(\NN_\eps(x))\mid x=y] .
\end{eqnarray*}
Note that $g(y)+g(\sigma_j(y))=\tilde g(y)+\tilde g(\sigma_j y)$,
but 
$|g(y)-g(\sigma_j(y))|\leq |\tilde g(y)-\tilde g(\sigma_j y)|$.
This implies
$g(y)^2+g(\sigma_j(y))^2\leq \tilde g(y)^2+\tilde g(\sigma_j y)^2$.
By summing over $y$, we obtain $\E (g^2) \leq \E(\tilde g^2)$.
Since $\E g = \E\tilde g$, the last claim of the lemma now follows.
\ep
\medskip

\proofof{main}
Let $\A\subset\Omega_n$.
Set $g = \kappa_1\kappa_2\dots\kappa_n \chi_\A$. 
Then by Lemma~\ref{l.shift}, $g$ is monotone,
$\H(g)\leq\H(\A)$ and  for each $\eps>0$
we have $ \var(g,\eps) \geq \var(\A,\eps)$.
Moreover, $g$ takes only the values $0$ and $1$.
By applying Theorem~\ref{t.tala} for $g$, and using 
Theorem~\ref{p.varphi}, Theorem~\ref{main} immediately follows.
\ep
\medskip

\proofof{p.mconverse}
Observe that for a monotone $f:\Omega_n\to \R$
\begin{equation}
\label{e.monfour}
I_j(f)=2|\hat f(\{j\})|,
\end{equation}
and therefore
\begin{equation}
\label{e.monfourh}
\H(f)=4\sum_j \hat f(\{j\})^2.
\end{equation}
Hence~\ref{p.mconverse} follows from Theorem~\ref{t.foursens}.
\ep
\medskip

Note that~(\ref{e.monfourh}) implies the well-known inequality
\begin{equation}
\label{e.monh}
\H(\A)\leq 1
\end{equation}
for monotone events $\A$.

\begin{rem}
It is tempting to look for a simpler proof of Theorem~\ref{main},
along the following lines.  Using~(\ref{e.fjnrmp}) with $p=2$,
we find that
\begin{equation}\label{e.ii}
\H(f) =  \sum_{j=1}^n \left(\sum_{S\supset \{j\}} 4 \hat f(S)^2\right)^2
= 16 \sum_{S,S'} \hat f(S)^2 \hat f(S')^2 |S\cap S'|
,
\end{equation}
where $f=\chi_\A$ for some event $\A\subset\Omega_n$.
This expression is more complicated than~(\ref{e.monfourh}),
but is still valid when $\A$ is not monotone.
The fact that $f$ is the indicator function of an event is summarized
by the equation $f^2=f$.  In terms of the Fourier transform,
this translates to a convolution equation
\begin{equation}
\label{e.charf}
\hat f \ast \hat f = \hat f
.
\end{equation}
(By replacing $f$ with $2f-1$, this transforms to the simpler looking
$\hat f\ast\hat f= \chi_{\{\emptyset\}}$.)
One may suspect that there should be a direct argument
that uses only~(\ref{e.charf}) and~(\ref{e.ii}) to prove that
for every $k=1,2,\dots$
$$
\sum_{|S|=k}\hat f(S)^2 \to 0
$$
when $\H(f)\to 0$.  Then Theorem~\ref{main} would follow from
Theorem~\ref{t.foursens}.
\end{rem}

\section{Correlation with majority} \label{s.maj}

\subsection{Uniform weights}

Fix some $n\in\N$.
Recall the definition~(\ref{e.majdef}) of the majority function $\Maj_K$,
and set $\Maj=\Maj_n=\Maj_{\nn}$.

\begin{thm}
\label{t.majcor}
Let $f:\Omega_n\to [0,1]$ be monotone. Then
$$
I(f)\leq C \sqrt n\, \E( f\Maj)\,\left(1+\sqrt{-\log\E(f\Maj)}\right)
\,,
$$
where $C$ is some universal constant.
\end{thm}

\proof
Write $\br f(k)$ for the average of $f$ on the set
$\left\{x:\sum_jx_j=k\right\}$:
$$
\br f(k)= {n\choose k}^{-1}\sum_{|x|=k} f(x)
\,.
$$
Then 
\begin{equation}\label{e.majcor}
\E(f\Maj) =
2^{-n}\sum_{k>{n\over 2}}{n\choose k}
\left(\br f(k)-\br f(n-k)\right)
\,.
\end{equation}
Recall that $s_jx=(y_1,\dots y_n)$ where
$y_j=1-x_j$ and $y_i=x_i$ for $i\neq j$.
Then 
$$
I(f)=2^{-n}\sum_x\sum_j |f(x)-f(s_j x)|
.
$$
Since $f$ is monotone, $f(x)-f(s_jx)\geq 0$ when
$x_j=1$ and $f(x)-f(s_jx)\leq 0$ when $x_j=0$.
Hence the expression for $I(f)$ simplifies,
\begin{eqnarray}
I(f)
&=&
 2^{-n}\sum_x f(x) \Bigl(2|x|-n\Bigr)
\nonumber
\\
&=&
 2^{-n}\sum_k{n\choose k} \br f(k) (2k-n)
\nonumber
\\
&=&
 2^{-n}\sum_{k>{n\over 2}}{n\choose k} \left(\br f(k)-\br f(n-k)\right)
(2k-n)
.
\label{e.inf}
\end{eqnarray}
For any $\lambda\geq 0$
write $k(\lambda)=(n+\lambda\sqrt n)/2$.
Since $0\leq \br f(k)\leq 1$,
by comparing (\ref{e.inf}) and (\ref{e.majcor}),
we obtain the following estimate.
\begin{eqnarray}
I(f)
&\leq& (2 k(\lambda) -n)\,\E(f\Maj)
+ 2^{-n}\sum_{k>k(\lambda)}{n\choose k}\left(\br f(k)-\br f(n-k)\right)
(2k-n)
\nonumber
\\ & \leq &
\lambda\sqrt n\, \E(f\Maj)
+ 2^{-n}\sum_{k>k(\lambda)}{n\choose k} (2k-n)
\label{e.i}
.
\end{eqnarray}
Because there are constants $C_1,C_2>0$ such that
\begin{equation}
\label{e.need}
2^{-n}{n\choose k}(2k-n) \leq C_1 \exp\left(-{(2k-n)^2\over C_2n}\right)
\end{equation}
holds for every $n$ and $k$,
by choosing $\lambda= C_3\sqrt{-\log \E(f\Maj)}$,
where $C_3$ is a sufficiently large constant,
we get
$$
2^{-n}\sum_{k>k(\lambda)}{n\choose k} (2k-n)
\leq C_4 \sqrt n \, \E(f\Maj),
$$
and the theorem follows from~(\ref{e.i}).
\ep \medskip

Given a set $K\subset\nn$, let
$\Maj_K$ denote the majority function on the set $K$;
that is,
$$
\Maj_K(x)=
\left\{
\begin{array}{rl}
{}-1 & \mbox{if $\sum_{j\in K} x_j < |K|/2$}\,;\\
0 & \mbox{if $\sum_{j\in K} x_j = |K|/2$}\,;\\
1 & \mbox{if $\sum_{j\in K} x_j > |K|/2$}\,,
\end{array}\right.
$$
Also set,
$$
I_K(f)=\sum_{k\in K} I_k(f)
.
$$

\begin{cor}
\label{c.submaj}
Let $K\subset\nn$ and suppose that
$f:\Omega_n\to [0,1]$ is monotone.
Then
$$
I_K(f)\leq
C \sqrt{|K|}\, \E(f\Maj_K)\,\left(1+\sqrt{-\log\E(f\Maj_K)}\right)
\,,
$$
where $C$ is some universal constant.
\end{cor}

\proof
Set $m=|K|$, and assume, that $K=\{1,\dots,m\}$.
Given $z\in \Omega_m$, set
$$
f_K(z)=2^{m-n}\sum_{y\in\Omega_{n-m}} f(z,y).
$$
Then $f_K$ is monotone and $I(f_K)=I_K(f)$.
Consequently, the corollary follows from Theorem~\ref{t.majcor}.
\ep \medskip

\proofof{t.lambii}
Assume, with no loss of generality, that
\begin{equation}\label{e.dec}
I_{j+1}(f)\leq I_j(f)
\end{equation}
for all $j\in\{1,\dots,n-1\}$.
Cor.~\ref{c.submaj} implies that
\begin{equation}\label{e.sfk}
\sum_{j=1}^k I_j(f)
\leq C_1 \Lambda(f) \left(1+\sqrt{-\log \Lambda(f)}\right)\sqrt k
\end{equation}
for some constant $C_1$ and every $k\in\nn$.
Subject to these constraints and~(\ref{e.dec}),
$\H(f)$ is maximized if equality occurs in~(\ref{e.sfk})
for every $k$.  Therefore,
\begin{eqnarray*}
\H(f)
& \leq &
C_1^2 \Lambda(f)^2 \left(1+\sqrt{-\log \Lambda(f)}\right)^2
\sum_{k=1}^n \left(\sqrt k-\sqrt{k-1}\right)^2
\\
&=&
O(1)
\Lambda(f)^2 \Bigl(1-\log \Lambda(f)\Bigr)
\sum_{k=1}^n k^{-1}
\\ & = &
 O(1) \Lambda(f)^2 \Bigl(1-\log \Lambda(f)\Bigr) \log n
.
\end{eqnarray*}
This proves the first part of Theorem~\ref{t.lambii}.
The second part now follows from Theorem~\ref{main}.
\ep
\medskip

Theorem~\ref{t.lambii} tells us that if $\Lambda(\A_m)\to 0$
fast enough for monotone events $\A_m$, then they are
asymptotically noise sensitive. 
Conversely, if a sequence of (not necessarily
monotone) events satisfies $\inf_m\Lambda(\A_m)>0$,
then it is not asymptotically noise sensitive.
This can be proven directly, and also follows from
Lemma~\ref{l.corstab} below.

It is interesting to note that

\begin{thm}
\label{t.ccfs}
Majority maximizes $I$ among monotone events $\A\subset\Omega_n$.
\end{thm}

This follows from~\cite{CCFS}, although the explicit statement
does not appear there. It also follows from the classical 
Kruskal-Katona theorem. See also \cite[Lem.~6.1]{FK}.

\subsection {General weights}

We will investigate now some relations between noise-sensitivity and 
weighted majority functions. 
Several of the properties we need for weighted
majority functions are easy to establish if the distribution of weights
allows us to use 
a normal approximation for $f(x)=\sum_j w_j x_j$. But, as it turns out,
working with arbitrary weights is harder.

Our first goal is to show that weighted majority functions are uniformly 
noise stable. This will imply the  ``only if" part of 
Theorem  \ref{t.lambii2}.  For this, the following easy (and quite
standard) lemma will be needed.

\begin{lemma}
\label {l.ll}
Let $w=(w_1,\dots,w_n)\neq 0$ and $f(x)=\sum_j w_j(2 x_j-1)$.
Then 
\begin{equation}
\label{e.p1}
\P[|f| \geq t \|w\|_2] \le 3t^{-4} 
\,,
\end{equation}
and
\begin{equation}
\label{e.p2}
\P [ |f| \leq 0.3\|w\|_2 ] \le 0.92
\,.
\end{equation}
\end{lemma}

A much stronger estimate than (\ref{e.p1}) is known
(see \cite{Petrov}).

\medskip
\proof
Without loss of generality, we assume that  $\|w\|_2=1$.
Then $\E[f^4]=3\|w\|_2^4-2\|w\|_4^4\leq 3$.  Hence (\ref{e.p1}) follows:
$$
\P[|f|\geq t] = \P[f^4\geq t^4] \leq t^{-4} \E[f^4] = 3 t^{-4}.
$$
This implies
$$
\E[{\mathbf1}_{\{f^2>t\}}f^2 ]
= \P[f^2>t] + \int_{s\geq t} \P[f^2>s]\,ds
\leq 3 t^{-2} + t^{-3}
\,.
$$
Hence 
$$
\E[{\mathbf 1}_{\{f^2\leq t\}}f^2 ]
=
\E[f^2]- \E[f^2 {\mathbf 1}_{\{f^2>t\}}]
\geq 1- 3 t^{-2} - t^{-3}
\,.
$$ 
We choose $t=10$, and obtain
\begin{multline*}
10 - 9.9\P[f^2\leq 1/10] 
=
10\P[f^2>1/10]+ \P[f^2\leq 1/10]/10
\\
\geq
\E[{\mathbf 1}_{\{f^2\leq 10\}}f^2 ] \geq 9/10
\,,
\end{multline*}
which gives (\ref{e.p2}).
\ep 
\medskip

\begin{lemma}
\label{l.sums}
Let $b>0$,
let $v_1,\dots,v_d\geq b$, and let $g= \sum_{j=1}^d z_j v_j$,
where $\P[z_j=1]=\P[z_j=-1]=1/2$, and the $z_j$ are independent.
Then for every $t\geq 1$ and every $s\in\R$, 
\begin{equation}
\label{e.sums}
\P\bigl[|g-s|\leq t b\bigr] \leq c\cdot t/\sqrt{d}
\,,
\end{equation}
where $c$ is some universal constant.
\end{lemma}

This lemma is a consequence of Theorem 2.14 in \cite{Petrov}, for example.
However, since the proof of that theorem is arduous, we now present a simple
combinatorial proof.

\medskip \proof
Let $x$ be a random uniform element in $\Omega_d$,
and let $\pi$ be a random uniform permutation of $\{1,2,\dots,d\}$.
Let $C$ be the collection of sets $S$ that have the form $S=\{j:\pi(j)<r\}$
for some $r\in\R$.  Then there is a unique $y\in C$ with $|y|=|x|$.
 Observe that $y$ is a random uniform element of $\Omega_d$.
Consequently, the distribution of $g$ is the same
as the distribution of $h(y):=\sum_{j=1}^d (1-2 y_j) v_j$,
where $y_j$ is $1$ or $0$ when $j\in y$ or $j\notin y$, respectively.
Since $C$ is totally ordered by inclusion,
there is at most one $S\in C$ such that $|h(S)-s|< b/2$.
So when $\pi$ is fixed, the probability that $|h(S)-s|\leq b/3$
is at most $\max\Bigl\{\P\bigl[|x|=r\bigr]:r\in \R\Bigr\}= O(1)/\sqrt{d}$.
This establishes (\ref{e.sums}) for $t=1/3$.  The result for
general $t\geq 1$ follows by applying the result for $t=1/3$ for
an appropriate succession of values of $s$.
\ep
\medskip

\proofof{t.stab}
Let $w=(w_1,\dots,w_n)\neq 0$ and $s_0\in\R$.
Let $f(x):=\sum_{j=1}^n w_j(2 x_j -1)$, and
consider the event $\Cal M :=\bigl\{x\in\Omega:f(x)>s_0\bigr\}$.
Take $\eps>0$, and
let $J\subset\nn$ be a random subset, where each $j\in\nn$
is in $J$ with probability $\eps$, independently.
Set $Y(J):= \sum_{j\in J} w_j(2 x_j -1)$.
Then $2Y(J)$ has the distribution of $f-\NN_\eps f$.
Let $\delta\in(0,1)$ and set
\begin{equation}
\label{e.adef}
a:=\inf\Bigl\{t>0:\P\bigl[\lvert Y(J)\rvert\geq t\bigr]\leq\delta\Bigr\}.
\end{equation}

Our goal is to give an estimate from above to
$P\bigl[|f|< 2a\bigr]$ in terms of $\eps$ and $\delta$,
which will tend to zero when 
$\delta$ is positive and fixed and $\eps\to 0$.

Set $W(J):= \sum_{j\in J} w_j^2$.  This is
the variance of $Y(J)$ conditioned on $J$.
Note that
$$
\P\bigl[\lvert Y(J)\rvert\geq a\mid J\bigr] = \P\bigl[Y(J)^2>a^2\mid J\bigr]
\leq
\E\bigl[Y(J)^2\mid J\bigr]/a^2 = W(J)/a^2.
$$
Therefore, 
\begin{multline}
\delta=\P\bigl[\lvert Y(J)\rvert>a\bigr] =
\sum_{X\subset \nn} \P\bigl[\lvert Y(J)\rvert>a\mid J=X\bigr]\P[J=X]
\\
\leq\sum_{X\subset\nn} \min\bigl\{1,a^{-2}W(X)\bigr\}\P[J=X]
= \E\Bigl[\min\bigl\{1,a^{-2}W(J)\bigr\}\Bigr]
\,,
\end{multline}
and we conclude that
\begin{equation}
\label{e.w}
\P\bigl[W(J)\geq\delta a^2/2\bigr]\geq \delta/2
\,.
\end{equation}

Now let $z_1,z_2,\dots,z_n$ be independent variables
that are uniform in $[0,1]$, and are independent from
$(x_1,\dots,x_n)$.  Let $m$ be the largest integer
such that $m\eps\leq 1$. Let $I_1,\dots,I_m$ be
disjoint open intervals in $[0,1]$, each of
length $\eps$.  Let $I_0:=[0,1]-\cup_{k=1}^m I_k$.
Let $J_k$ ($k=0,1,\dots,m$)
be the set of $i\in\nn$ with $z_i\in I_k$.
Then each $J_k$ with $k>0$ has the same distribution as
$J$ above.  Let $\A_k$ be the event
that $W(J_k)\geq \delta a^2/2$.  Then from (\ref{e.w})
with $J_k$ in place of $J$ we find that $\P[\A_k]\geq\delta/2$
for $k=1,2,\dots,m$.

We claim that for $k\neq k'$ the events $\A_k$ and $\A_{k'}$ are
negatively correlated.  This can be established by
proving by induction on $n$ that the events $W(J_k)\geq s_1$
and $W(J_{k'})\geq s_2$ are negatively correlated for each $s_1,s_2\in\R$
(which is intuitively obvious, since the intervals $I_k$ and
$I_{k'}$ are disjoint).
Let $K$ be the number of $k>0$ such that the
event $\A_k$ occurs.
Then
$$
\E[K] = \sum_{k=1}^m \P[\A_k] \geq m\delta/2
\,,
$$
and
$$
\E[K^2]-\E[K]^2 =\sum_{k,k'}\P[\A_k\cap\A_{k'}]-\left(\sum_k
\P[\A_k]\right)^2
\leq \E[K],
$$
because the events $\A_k,\A_{k'}$ are negatively correlated
when $k\neq k'$.
Therefore
\begin{equation}
\label{e.Kest}
\begin{split}
\P[K< m\delta/4]
&
\leq  \P\left[(K-\E K)^2 > (m\delta/4)^2\right]
\\&
\leq 4 (m\delta)^{-2} \E\left[(K-\E K)^2\right] 
= 4 (m\delta)^{-2} \left(\E[K^2] -(\E K)^2\right)
\\&
\leq 4 (m\delta)^{-2} \E[K]
\leq 4 m^{-1}\delta^{-2}
\,.
\end{split} 
\end{equation}
Let $L$ be the set of $k\in[m]$ such that $|Y(J_k)|> a\sqrt{\delta}/10$.
By (\ref{e.p2}), applied to $Y(J_k)$ in place of $f$,
$$
\P[k\in L\mid\A_k]\geq 8/100
.
$$
Moreover, conditioned on all the $J_k$, the events $\{k\in L\}$ are
independent.  Consequently, a calculation similar to (\ref{e.Kest})
gives
$$
\P\bigl[|L|<  m\delta/100\mid K\geq m\delta/4\bigr] \leq 
O(1)m^{-1}\delta^{-2}
.
$$
When we use this and (\ref{e.Kest}) together, we get
\begin{equation}
\label{e.lbb}
\P\bigl[|L|<  m\delta/100\bigr] \leq 
O(1)m ^{-1} \delta^{-2} .
\end{equation}
If we condition on $L$, on all $Y(J_k)$ for $k\notin L$ and
on all $|Y(J_k)|$ for $k\in L$, then what remains to determine
$f$ are only the signs of $Y(J_k)$ with $k\in L$.  Moreover,
these signs are independent, and are $+$ or $-$ with
probability $1/2$.
Hence we may apply Lemma~\ref{l.sums} with 
$b :=  a\sqrt{\delta}/10$, $d:=|L|$,
$s:= s_0-\sum_{k\notin L} Y(J_k)$,
$g=\sum_{k\in L} Y(J_k)$, and take
$v=(v_k)$ to be the sequence $\bigl(|Y(J_k)|:k\in L\bigr)$.
The conclusion is that for $t\geq 1$
$$
\P\Bigl[|f-s_0|\leq  ta\sqrt{\delta}/10\ \Big|\ 
|L| \geq  m \delta/100\Bigr]
\leq  O(1) t / \sqrt{ m\delta } 
\,.
$$
Together with (\ref{e.lbb}), (and choosing $t=20/\sqrt{\delta}$)
this gives
\begin{equation}
\P\Bigl[|f-s_0|\leq 2a\Bigr] 
\leq O(1)\left(\eps\delta^{-2}+\sqrt{\eps}\delta ^{-1} \right).
\end{equation} 

We now come to analyze the effect of noise.
Because $2Y(J)$ has the same distribution as
$f-\NN_\eps f$, for every $a>0$
$$
\P\bigl[\Cal M\triangle \NN_\eps\Cal M\bigr]
\leq \P\bigl[\lvert f-s_0\rvert \leq 2a\bigr]
+\P\bigl[\lvert Y(J)\rvert\geq a\bigr]
$$
Choose $\delta:=\eps^{1/4}$ and, as before, use (\ref{e.adef})
to define $a$.
Then $\P\bigl[|f-s_0|\leq 2a\bigr]\leq O(1)\eps^{1/4}$
and $\P\bigl[|Y(J)|\geq a\bigr]\leq \eps^{1/4}$.
Consequently,
\begin{equation}
\label{e.miq}
\P\bigl[\Cal M\triangle \NN_\eps\Cal M\bigr]
\leq O(1)\eps^{1/4}
,
\end{equation}
and the theorem immediately follows.
\ep

\begin{que}
What is the best exponent possible on the right hand side of~(\ref{e.miq})?
\end{que}

\begin {rem2}
Yuval Peres and Elchanan Mossel found a simple proof showing that, 
as expected, the correct exponent
is 1/2.
\end {rem2}

\begin{rem}
It follows from Theorems~\ref{t.stab} and~\ref{main} that
$\inf \bigl\{\H(\Cal M) : \Cal M \in \mathfrak M\bigr\}>0$. (A direct 
proof will follow.)
We conjecture that $\H(\Cal M)$ is minimized among
$\Cal M_{w,0}\subset\Omega_n$ in $\mathfrak M$ when
all the weights are equal. 
It is a consequence of Theorems~\ref{t.stab} and~\ref{t.foursens}
that
$$
\lim_{k\to\infty}
\sup_{\Cal M\in\mathfrak M}
\sum_{|S|>k}\hat \chi_{\Cal M}(S)^2 =0.
$$
We actually expect that among weighted majority events in
$\Omega_n$, the one with equal weights is the least stable,
and for every $k>1$ maximizes  $\sum_{|S|>k}\hat \chi_{\Cal M}(S)^2$.
\end{rem}

For the proof of \ref{t.lambii2}, the following will be needed.

\begin{lemma}
\label{l.corstab}
Let $\A_m,\Cal B_m\subset\Omega_{n_m}$ be two sequences of
events.  Suppose that the sequence $\{\A_m\}$ is noise-sensitive,
while the sequence $\{\Cal B_m\}$ is noise-stable.
Then 
$$
\lim_m\P[\A_m\cap\Cal B_m] -\P[\A_m]\P[\Cal B_m]= 0 .
$$
\end{lemma}

\proof
This can be proven directly, but
since $\P[\A_m\cap\Cal B_m] = \E\bigl[\chi_{\A_m}\chi_{\Cal B_m}\bigr]$,
the lemma is immediate from~\ref{t.foursens}.
\ep
\medskip

Let the {\bf influence vector} of an event $\A\subset\Omega_n$ be the
vector $I^{\A}:=\bigl(I_1(\A),\dots,I_n(\A)\bigr)\in\R^n$.

\medskip
\proofof {t.lambii2} 
The ``only if'' direction follows from Theorem~\ref{t.stab}
and Lemma~\ref{l.corstab}.

For the  other direction, we need to show that  
monotone, noise-insensitive events $\A\subset\Omega_n$
have a non-vanishing correlation with some 
weighted majority event $\Cal M_w$, $w\in [0,1]^n$.
Talagrand's Theorem 1.1 \cite {T2} gives a lower bound on
the correlation of monotone events.
This theorem asserts, in particular, that for two monotone events,
if the inner product of their influence vectors is bounded away 
from zero, then the correlation between them is also bounded away from zero%
\footnote {For uniformly stable events, it seems that also the 
converse is true: if the correlation is bounded away from zero,
then so is the inner product of their influence vectors. For monotone
uniformly stable events, this follows from the two-event 
version of Theorem~\ref {t.tala}.}.

We know from Theorem~\ref{main} that for noise-insensitive events,
$\bigl\|I^{\A}\bigr\|_2$ is bounded away from zero.
It remains to show that for every $v\in[0,1]^n$ with $\|v\|_2=1$,
we can find a weighted majority function 
$\Cal M=\Cal M_{w}$, $w\in[0,1]^n$, such that the inner product
$\left<I^{\Cal M},v\right>$
is bounded away from zero. We will prove that this holds when one
chooses $w:=v$.
\medskip

Given any $w\in\R^n$, $w\neq 0$, let $I^w\in\R^n$ denote the
influence vector of $\Cal M_{w}$, $I^w_j := I_j(\Cal M_{w})$.

\begin{prop}
\label{p.infvec}
There is an absolute constant $c>0$ such that
$\left< w , I^w\right>\geq c$
for every $n=1,2,\dots$ and every $w\in\R^n$
with nonnegative coordinates and $\|w\|_2=1$.
\end{prop}

\proof
Set $f(x)=\sum_{j=1}^n (2x_j-1)w_j$ for $x\in\Omega_n$.
Then $\hat f(\{j\})=w_j$ for $j\in\nn$ and $\hat f(S)=0$ for
$S\subset\nn$, $|S|\neq 1$.
On the other hand, $I_j(\Cal M_w)=\hat{M}_w(\{j\})$,
where $M_w=\hbox{\rm sign}(f)$.
Therefore,
$$\left< w , I^w\right>=
\bigl<\hat f , \hat{M}_w \bigr>=
\left<f , {M}_w \right>=\E\left|f(x)\right|
\,,
$$
which is bounded from below, by (\ref{e.p2}).
This completes the proof of the proposition,
and the proof of Theorem~\ref{t.lambii2}.
\ep\medskip

\begin{rem}\label{r.nouni}
We now show that one cannot remove the $\log$ in
Theorem~\ref{t.lambii}. Fix some  $k,n\in\Z$ with $n\geq k>0$.
Let $w_j=1/\sqrt {j\log n}$ for $j=1,\dots,n$,
and let $u_j=1/\sqrt k$ for $j\leq k$ and $u_j=0$ for $j> k$.
Set $f_w(x)=\sum_{j=1}^n (2x_j-1)w_j$ and
$f_u(x)=\sum_{j=1}^n (2x_j-1)u_j$, where $x\in\Omega_n$.
Then the event $f_w\geq 0$ is noise stable, by~\ref{t.stab}.
We show that $\P[f_w\geq 0\mid f_u\geq 0]\to 1/2$ as
$n\to\infty$, no matter how $k=k(n)$ is chosen.
Indeed, given any $x\in\Omega$, let $s(x):=\sqrt k f_u(x)=
\sum_{j\leq k}(2 x_j-1)$.  If $s(x)<0$, let
$\overline x$ be obtained from $x$ by replacing $-s(x)$
of the $0$ entries in $x$ by $1$'s, where the set of entries replaced
is chosen randomly and uniformly among all possibilities, and if $s(x)\geq
0$,
set $\overline x=x$.
Then $\P[f_w(x)\geq 0\mid f_u(x)\geq 0] =\P[f_w(\overline x)\geq 0]$.
Therefore, by Lemma~\ref{l.sums} applied to $w$, it is enough to show that
$f_w(\overline x)-f_w(x)\to 0$ in probability as $n\to\infty$.
This follows from 
$$
\E\bigl\lvert f_w(\overline x)-f_w(x)\bigr\rvert  =
(2/k)\E\bigl[\max\{0,-s(x)\}\bigr]\sum_{j=1}^k w_j
= O(1)/\sqrt{\log n}.
$$
\end{rem}
\bigskip

\section {An application to percolation} \label{s.perc}

Let $R$ be an $(m+1)\times m$ rectangle in the square grid $\Z^2$,
and let $\Omega$ be the set of all functions
from $E$, the set of edges of $R$, to $\{0,1\}$.
We identify $\Omega$ with $\Omega_n$; where $n=n_m=|E|=2m^2-1$.
A point  $x\in\Omega$ is called a {\bf configuration},
and can be identified with the subgraph consisting
of all vertices of $R$ and all edges $e$ with $x(e)=1$.
A connected component of this graph is called a {\bf percolation cluster}.

Let $\cro=\cro_m\subset\Omega$ be the event that there is
a left-right crossing of $R$; that is, $\cro$
is the set of all configurations that contain a path
joining the left and right boundaries of $R$.
An easy and well known application of duality shows that $\P[\cro]=1/2$.

Kesten~\cite{Ke} gives an estimate from above for the probability
that an edge near the middle of $R$ is pivotal for $\cro$.
Similar estimates for edges near the boundary can probably be
extracted from Kesten's paper.  These give an inequality
of the form $I_j(\cro_m)\leq m^{-1-c}$, $c>0$, for each $j$.
Then Theorem~\ref{t.weak} implies~\ref{t.perc}.
However, we prefer to present another proof,
based on Theorem~\ref{t.lambii}.

The only percolation background needed to understand the proof is
that in our situation the probability that a vertex in $R$ is
connected in the configuration to some vertex at Euclidean
distance $r$ is at most $Cr^{-1/\rho}$, for some constants $C,\rho>0$.
This follows from the celebrated Russo-Seymour-Welsh Theorem~\cite{R,SW}
(see also~\cite{G}).

\proofof{t.perc}
Let $E_r$ be the set of edges in the right half of $R$,
with edges exactly centered included.  Let $K\subset E_r$.
We now estimate $\E(\chi_\cro\Maj_K)$.

\def\vis{{\rm VISITED}}

Consider the following algorithmic method of randomly selecting
a configuration.  Let $\omega^K$ and $\hat\omega^{K}$ be two
independent elements of $\Omega_{|K|}$ and $\Omega_{n-|K|}$,
respectively.  Let $V_1$ be the set of vertices on the left
boundary of $R$, and set $\vis=\emptyset$.
As long as there is some edge $[v,u]\notin\vis$
joining a vertex $v\in V_1$ to a vertex $u\notin V_1$, choose some such
edge $e=[v,u]$, and do the following. 
Append $e$ to $\vis$.
If $e\in K$, let $y(e)$ be the first bit in the sequence
$\omega^K$ that has not been previously used by the algorithm,
while if $e\notin K$ let $y(e)$ be the first bit in the
sequence $\hat\omega^K$ that has not been previously used by
the algorithm.  If $y(e)=1$,
then adjoin to $V_1$ the vertex $u$.

This procedure defines $y$ for all $e\in\vis$.
Let $z\in\Omega$ be random, uniform, and independent of $y$,
and let $x=y$ on $\vis$ while $x=z$ on $E-\vis$.
This defines a configuration $x\in\Omega$.

The following is obvious:

\begin{lemma}\label{l.obv}
The configuration $x$ given by the above algorithm is
uniformly distributed in $\Omega$.  The event
$x\in\cro$ is equal to the event that at the end
of the algorithm $V_1$ intersects the right boundary
and is independent from $z$ (can be determined by $y$).
\ep
\end{lemma}

Let us estimate the probability that $K\cap\vis$ is large.
An edge $e\in K$ is in $\vis$ iff there is in $x$ a path
joining a vertex of $e$ to the left boundary of $R$.
Since $K\subset E_r$,
it follows from the above stated consequence of the Russo-Seymour-Welsh
Theorem that the probability for the latter event is 
bounded by $C m^{-1/\rho}$, for some constants $C,\rho>0$.
Consequently, 
$$
\E|K\cap\vis| \leq C |K| m^{-1/\rho},
$$
which implies
$$
\P[\A_1]\leq C m^{-1/(3\rho)}
,
$$
where $\A_1$ is the event
$$
\A_1 := \Bigl\{x\in\Omega: |K\cap\vis|\geq |K|m^{-2/(3\rho)}\Bigr\}
.
$$
Let $\A_2$ be the event that there is an integer
$j$ in the range $1\leq j\leq |K| m^{-2/(3\rho)}$
such that 
$$
\left|{j\over2}-\sum_{i=1}^j \omega^K_i\right| \geq
\sqrt{|K| m^{-2/(3\rho)}} \log m.
$$
It is easy to see that the $\P[\A_2]$ decays super-polynomially in
$m$; in particular,
$$
\P[\A_2]\leq O( m^{-1/\rho} )
.
$$
As $\P[\A_1\cup\A_2]\leq O(1) m^{-1/(3\rho)}$, we have 
\begin{equation}\label{badset}
\E(\chi_{\A_1\cup\A_2} \chi_\cro\Maj_K) \leq O(1) m^{-1/(3\rho)}
.
\end{equation}
Now suppose that the algorithm produced a $y$
such that $\A_1\cup\A_2$ does not hold.
Then it follows that
$$
\left|{|\vis\cap K|\over 2}-\sum_{e\in\vis\cap K} y(e)\right|
\leq 
O(1)
\sqrt{|K| m^{-2/(3\rho)}} \log m.
$$
This implies that 
$$
\E[\Maj_K(x)\mid y]\leq
O(1) m^{-1/(3\rho)}{\log m}
,\qquad\forall y\notin\A_1\cup\A_2
.
$$
Since $x\in\cro$ can be determined from $y$, we get
$$
\E\Bigl((1-\chi_{\A_1\cup\A_2})\chi_\cro \Maj_K \Bigr)
\leq
O(1)m^{-1/(3\rho)}{\log m}
.
$$
In view of~(\ref{badset}) this implies
$$
\E(\chi_\cro \Maj_K) \leq O(1)m^{-1/(3\rho)}{\log m}
,
$$
and Cor.~\ref{c.submaj} gives
\begin{equation}
\label{e.IK}
I_K(\cro)\leq O(1)\sqrt{|K|} m^{-1/(3\rho)}({\log m})^{3/2}
\end{equation}
for every $K\subset E_r$, since $\cro$ is monotone.
By symmetry, this would also hold for $K\subset E-E_r$,
and therefore for every $K\subset E$.
Consequently, by the proof of Theorem~\ref{t.lambii}
\begin{equation}\label{e.hperc}
\H(\cro)\leq O(1) m^{-2/(3\rho)}(\log m)^4
.
\end{equation}
An appeal to Theorem~\ref{t.weak} completes the proof.
\ep

\begin{rem}
Since $I(\cro)=\sum_eI_e(\cro)$
is also the expected number of pivotal edges
for $\cro$, (\ref{e.IK}) shows that
the expected number of pivotal edges is bounded by
$$
O(1)m^{1-1/(3\rho)}(\log m)^{3/2}
.
$$
Although this is better than the general bound of
$O(1)m$ that follows from Theorem~\ref{t.ccfs},
a somewhat better bound can be extracted from Kesten's~\cite{Ke}.
\end{rem}

\begin{cor}
There is a constant $c>0$ with the following property.
If $\eps=c/\log m$, then for large $m$, with probability at least $1/4$,
$\Bigl|\{x,\NN_\eps(x)\}\cap \cro\Bigr|=1$.  That is, if each edge
is switched with probability $c/\log m$, independently,
then the crossing is likely to be created or destroyed.
\end{cor}

The corollary follows from~(\ref{e.hperc}) and Theorem~\ref{t.pow}.
The details are left to the reader.

\section {Some conjectures and problems concerning percolation}
\label{s.percconj}

\subsection{Other sensitivity conjectures} 

Consider the crossing event $\cro_m$ for an $(m+1) \times m$ 
rectangle in the square grid $\Z^d$.
By Theorem~\ref{t.perc} and Section~\ref{s.noise}, from knowing which edges
are open for 
all but a small random set of edges, we have almost no information whether 
crossing occurs. This suggests that for some deterministic
subsets of the rectangle $R=R_m$, knowing the configuration restricted to
that configuration typically gives almost no information whether crossing
occurs. It follows from the Russo-Seymour-Welsh Theorem~\cite{R,SW} that
 $E_r$, the set of edges in the right half of the rectangle,
 is not such a subset. Yet we believe that all
the horizontal edges (or all the vertical edges) is such a subset. 
That is, let $x,y \in \Omega$ be two independent uniform-random
configurations.   Let $z(e)=x(e)$ for horizontal edges $e$, and $z(e)=y(e)$
for
vertical edges.
Let $p(\omega)=\P[z \in \C | x = \omega]$. 

\begin{conj}
\label{con.horiz}
For any $\eps>0$, for all sufficiently large $m$,
$$
\P \Bigl\{\omega \in \Omega : |p(\omega)-1/2| > \eps\Bigr\} < \eps.
$$  
\end{conj}
\medskip

Here is a variant of this conjecture for {\bf Voronoi percolation}. 
Fix a square in $\R^2$. Voronoi percolation is performed in two steps. 
First pick $n$ points in the square uniformly and independently. 
Second each cell in the Voronoi tessellation determined by the
chosen points is declared open with probability $1/2$, and closed
otherwise, independently of the other cells, 
(see Benjamini and Schramm~\cite{BS} for the exact definitions and a study
of 
Voronoi percolation). By duality, the probability of open left-right open 
crossing is $1/2$.  In the spirit of Theorem~\ref{t.perc}, we conjecture
that
typically, knowing the Voronoi tessellation (but not knowing which cells are
open) gives
almost no information whether an open left-right crossing exits.

\subsection {Stronger sensitivity conjectures}

As before, let $\cro_m$ denote the crossing event for an $(m+1) \times m$

rectangle in the square grid $\Z^2$. 
\begin {conj} 
\label {c:strong}
There exists a $\beta >0$ so that  
$\lim_{m\to\infty}\phi (\cro_m, m^{-\beta}) =0$. 
\end {conj}

It is known \cite {KZ,Ke}, respectively, that for some reals $0 < b_1 <b_2<1$, 
$$m^{b_1} \le I(\cro _m) \le m^{b_2}, $$ 
and it is conjectured (see, e.g., \cite {BH}, p.91)  
that $I(\cro _m)$ behaves like $m^{3/4}$. 

\begin {prob} Is it true that
$\lim_{m\to\infty}\phi (\cro_m, \epsilon_m)=0$ 
when $\epsilon_m=o(m^{-3/4})$? 
\end {prob}

Recall that our proof of noise sensitivity for $\cro_m$ used (indirectly)  
upper bounds on $I(\cro _m)$.
On the other hand, there is 
a simple heuristic argument that directly relates 
noise sensitivity to the {\it lower bounds} on $I(\cro _m)$ 
and the distribution of the number of 
pivotal edges. Namely, we cannot expect the crossing event 
$\cro_m$ to be stable under noise which, with a very high probability,
will flip pivotal edges. This argument tends to support 
Conjecture \ref {c:strong}.

\subsection {Dynamical percolation} 

Dynamical percolation was 
introduced by H\"{a}ggstr\"{o}m, Peres and Steif~\cite{hag96}.
Consider the following process.  
Let $\{X_e\}$ be independent Poisson point processes in $\R$ indexed by the
edges $e\in E_R$ of the $(m+1)\times m$ rectangle $R=R_m$ in $\Z^2$.
Let $x_0:E_R\to\{0,1\}$ be random-uniform.  For each $t>0$ set
$x_t(e):=x_0(e)$ if the number of points in $(0,t]\cap X_e$ is even,
and $x_t(e):=1-x_0(e)$ if the number is odd.
This gives a continuous time stationary Markov chain $x_t$ in
$\Omega=\{0,1\}^{E_R}$.
Write $\Pdy$ for the probability measure governing this process. 
For each fixed $t$, the random variable $x_t$ can be thought of
as ordinary (Bernoulli($1/2$)) percolation in $\Z^2$.

An interesting problem raised by~\cite{hag96} is weather
there are (exceptional, random) times $t$ in which there is
an infinite percolation cluster in $x_t$.  The result described
below might be relevant.

As before, let $\cro_m$ denote the set of configurations in
$\Omega$ that have an open left-right crossing of $R_m$.
For all $t$, $\Pdy[x_t\in\cro_m] =1/2$.   
Let $S_m$ be the set of switching times; that is,
$S_m$ is the boundary of $\{t\geq 0: x_t\in\cro_m\}$.
As a corollary of Theorem~\ref{t.perc}, we have,

\begin{cor}
\label{c.dyn}
$\Bigl|S_m\cap[0,1]\Bigr|\to\infty$ in probability.
\end{cor}

\proof
Suppose $s>t\geq 0$.  
Observe that the distribution of the pair $(x_t,x_s)$ is the same
as the distribution of the pair $\Bigl(x_0,\NN_\eps(x_0)\Bigr)$,
where $\eps$ is a function of $s-t$ and $\eps>0$ when $s>t$.
(Actually, $\eps/(s-t)\to 1$ as $s-t\to0$.)

Let $k$ be some positive integer, and set $\eps=\eps(1/k)$.
Let $t_j:=j/k$.  Let $\Cal W$ be the set of $\omega\in\Omega$ such that
$\Bigl|\P[\NN_\eps(\omega)\in \cro]-1/2\Bigr|>1/4$.
Then $\P[\Cal W]\to 0$ as $m\to\infty$, by Theorem~\ref{t.perc}.
Let $\Cal Z(a,b)$ be the event that $S\cap[a,b]=\emptyset$.
Observe that for $\omega\notin\Cal W$, we have
$$
\Pdy\Bigl[\Cal Z(t_j,t_{j+1})\mid x_{t_j}=\omega\Bigr]\leq 3/4
\,,
$$ 
because $\Cal Z(t_j,t_{j+1})$ is disjoint from the event
$\Bigl|\{x_{t_j},x_{t_{j+1}}\}\cap\cro\Bigr|=1$.
Hence we can make the following estimate,
\begin{eqnarray*}
&&\!\!\!\!\!\!\!\!\!\!\!\!\!\!\!\!\!\!\!\!\Pdy[\Cal Z(0,t_{j+1})]
=
\Pdy[\Cal Z(0,t_{j})\cap\Cal Z(t_j,t_{j+1})]
\cr &=&
\sum_{\omega\in\Omega}
\Pdy\Bigl[\Cal Z(0,t_{j})\cap\Cal Z(t_j,t_{j+1})\mid
   x_{t_j} =\omega\Bigr]\P\{\omega\}
\cr &=&
\sum_{\omega\in\Omega}
\Pdy\Bigl[\Cal Z(0,t_{j})\mid x_{t_j} =\omega\Bigr]
\Pdy\Bigl[\Cal Z(t_j,t_{j+1})\mid x_{t_j} =\omega\Bigr]\P\{\omega\}
\cr && \qquad\qquad\hbox{(by the Markov property for $x_t$)}
\cr &\leq &
 \P[\Cal W] +
\sum_{\omega\in\Omega-\Cal W}
\Pdy\Bigl[\Cal Z(0,t_{j})\mid x_{t_j} =\omega\Bigr]
\Pdy\Bigl[\Cal Z(t_j,t_{j+1})\mid x_{t_j} =\omega\Bigr]\P\{\omega\} 
\cr &\leq &
\P[\Cal W] +(3/4) \sum_{\omega\in\Omega-\Cal W}
\Pdy\Bigl[\Cal Z(0,t_{j})\mid x_{t_j} =\omega\Bigr]\P\{\omega\}
\cr &\leq &
\P[\Cal W] +(3/4) \sum_{\omega\in\Omega} 
\Pdy\Bigl[\Cal Z(0,t_{j})\mid x_{t_j} =\omega\Bigr]\P\{\omega\}
\cr &=&
 \P[\Cal W] +(3/4) \Pdy\Bigl[\Cal Z(0,t_j)\Bigr].
\end{eqnarray*}
Using this inequality and induction gives
$\Pdy\Bigl[\Cal Z(0,t_j)\Bigr]\leq 4\P[\Cal W] +(3/4)^j$.
By stationarity, for every $t\geq 0$,
the same estimate for the probability of
$\Cal Z(t,t+j/k)$ holds.  Since $k$ may be chosen
arbitrarily large, and $\P[\Cal W]\to 0$ as $m\to\infty$,
the corollary easily follows.
\ep

\subsection {Limits and conformal invariance}

The motivating questions behind this work were the conjecture
regarding the existence of the limit 
and the conformal invariance conjecture for two-dimensional percolation.

These conjectures say, roughly, that the crossing probabilities
inside a domain between two boundary arcs have a limit
as the mesh of the grid goes to zero, and the limit is
invariant under conformal transformations of the domain and
the boundary arcs. 
For more details, 
see Langlands, Pouliot and Saint-Aubin~\cite{LPS}.  

Consider a triple $\Cal G=\left<G,A,B\right>$, where $G=(V,E)$ is a  finite
planar graph with $m$ edges, and $A,B\subset V$.
Let $p_{\Cal G}$ be the probability that there
is an open crossing from $A$ to $B$ in a uniform-random configuration
$x\in\Omega=\{0,1\}^E$. 

Let $\Cal H= \left<H,A',B'\right>$ be a triple obtained from $G$ by 
the following operation: for every edge $e$ of $G$ delete $e$ with
probability
$(1-t)/2$ contract $e$ with probability $(1-t)/2$ and leave $e$ unchanged
with probability $t$, independently of the other edges.  
$\Cal H $ is a random variable 
which takes values in planar graphs
with two distinguished vertex sets. 

Of course, $E(p_{\cal H})=p_{\cal G}$, and noise sensitivity, when
it applies, 
asserts that the value of $p_{\cal H}$ is concentrated around the mean.
Noise sensitivity enables one to relate the crossing probabilities
of percolation on different graphs and  
we had hoped that it will be relevant to  conjecture regarding the 
existence of the limit conjecture.  
At present, however, such
applications
are beyond our reach as we do not have a good understanding of planar graphs
which are obtained by random deletions and contractions of the form 
described above when ${\cal G}$ is a rectangle in the square grid. 

To be more specific, 
suppose that we take $G$ to be the 
$m\times cm$ rectangle in $\Z^2$ ($c>0$ some fixed constant) and let
$A$ and $B$ be its left and right boundaries. 
It follows from Theorem~\ref{t.pow} and~(\ref{e.hperc}) 
that  $p_{\Cal H}-p_{\Cal G}\to 0$ 
in probability, provided that $t \cdot \log m\to\infty$. 
(Conjecture \ref {c:strong} would give it even when $t \cdot 
m^\beta \to \infty$
for some $\beta >0$.) 
It is conjectured that the crossing probability 
tends to a limit as $m$ tends to infinity and an
approach to this conjecture
would be to relate the distribution of 
such random planar graphs starting from similar rectangles of different
sizes.
(The values of $t$ should depend on the size but be large enough 
that noise sensitivity applies).

In a different direction, 
the random planar graphs $\cal H$ obtained when you start with the 
$(m+1) \times m$ grid and let $t = m^{-3/4}$ are of special interest
and might be related to models of random planar graphs 
in mathematical physics \cite {ADJ}.

\subsection {Fourier-Walsh coefficients of percolation}

It is a natural question to try to understand the Fourier-Walsh
coefficients of boolean functions given by percolation problems.
Consider (again) the event $\cro=\cro_m$ of a left-right crossing
of an $(m+1)\times m$
rectangle $R=R_m$ of the 
square grid, $\Z^2$.
Let $f_m:=\chi_{\cro_m}$. 
The Fourier coefficients of $f_m$ are indexed by subsets of $E_R$, the edges
in $R_m$. The values $\hat f^2$ can be regarded as a measure on the space of
subgraphs of $R_m$. 

\begin {prob}
Describe this measure!
\end {prob}

It follows from Theorem~\ref{t.pow} and our estimates for $\H(\cro_m)$,
that all but a negligible part 
of the $L^2$ weight of the Fourier coefficients
$\hat f(S)$, where $S$ is non-empty, is for $|S|> c \log m$. 
Conjecture \ref {c:strong} is equivalent to the assertion that, 
in fact, this is true for $|S| > m^{\beta}$ for some
$\beta >0$. Conjecture~\ref{con.horiz} is equivalent to the statement that 
for all but a negligible part of these Fourier coefficients, the number of 
vertical edges in $S$ tends to infinity with $m$.

\subsection {Other models of statistical mechanics}

It would be of interest to extend
the results of this paper as well as earlier results on influence
(\cite {KKL, FK}) to other models of statistical mechanics, such as 
the Ising and Potts models. 
Many of the results on influence and on noise sensitivity 
should be extendible to measures on 
$\Omega_n$ for which the coordinate variables are {\bf positively
associated}, 
namely, measures for which every two monotone real functions 
are positively correlated.

\section {Some further examples} \label {s.ex}

We will discuss now four examples, the first two were considered by Ben-Or
and
Linial \cite {BL}. 

\subsection {Tribes}

Consider $n$ boolean variables divided into $t$ tribes $T_1, T_2 \dots, T_t$

of size $s$ each, and
let $f$ be the boolean function which take the value 1 if for some
$j$, $1 \le j \le t$, all variables of $T_j$ equal 1.
If $s =  \log n -\log \log n +\log \log 2$, then 
$\P[f=1] \approx \frac 12$. 
Also note that $I_k(f) \sim \log n / n$ for every $k$. 
It is easy to show directly that
$f$ will be immune to $\epsilon$-noise when $\epsilon=o(1/\log n)$
and will be devastated by $\epsilon$-noise if $\epsilon \log n \to \infty$.

Thus, $J(f) \sim \log n /n$. 

\subsection {Recursive majority on the ternary tree}

Consider $n=3^t$ boolean variables which form the leaves of a rooted
ternary tree of height $t$.
A boolean function $f$ is defined as follows: Given values for the
variable on the leaves compute for each other vertex its value as the 
majority of the values of its sons and set the value of $f$ to be the
value of the root.   

Ben Or and Linial showed that $I_k(f) \sim n^{- \log2/\log 3} $ for
every $k$ and thus $\alpha (f) \to 1- \log 2/ \log 3$ as
$t\to\infty$. It is easy to see that 
also $\beta (f) \to 1-\log 2 /\log 3$. This follows at once from the
following 
observation: for $t=1$, if we switch the 
value of each leaf with probability $p$
independently,
then for small $p$ the probability that the outcome will be switched is
$(3/2)p+o(p)$. 

\begin {conj}
There is an absolute constant $\beta_0 < 1/2$ (find it!){} such that 
for every monotone Boolean function $f$, $\beta (f) \le \beta_0$.
\end {conj}

\begin {rem2}
It was pointed out by Mossel and Peres 
by considering certain recursive majorities on larger trees 
that this conjecture is false. 
\end {rem2}

\subsection {Number of runs}

We considered mainly monotone events. Here is an interesting noise stable
non-monotone event. Given a string of $n$ bits $x_1, x_2, \dots, x_n$, 
let $R(x_1, x_2, \dots x_n)$ be the  
number of runs. 
Thus, $R$ is one plus the number of pairs
of consecutive variables with different values. 

The event that $R(x_1, x_2, \dots x_n)$ 
is larger than its median is noise stable. 
Indeed, write $y_i= x_i~ \oplus~ x_{i+1}$,  $~i=1, \dots ,n-1$, 
and note that the $y_i$'s are independent, and $R$ is 
just the majority event on the $y_i's$. 
( Here $\oplus$ is addition mod $2$; that is, xor.)

\subsection {Majority of triangles}

We considered only the case where $p$ is a constant. 
When $p$ tends to zero with $n$, new phenomena occur.
Consider, for example, random graphs on $n$ vertices 
with edge probability $p = n^{-a}$, $a>0$ and the event that the 
number of triangles in the graph is larger than  
its median. This is a noise stable event but its correlation with  
majority (or any weighted majority) tends to 0 as $n$ tends to infinity.

\section {Relations with complexity theory} \label{s.compl}

There is an interesting connection between the complexity of 
boolean functions and the notions studied in this paper.

\subsection {AC0 and influences}

An important complexity class AC0 of Boolean functions   
are those which can be expressed by Boolean circuits
of polynomial size (in the number of variables) and bounded depth.
Boppana \cite {Bop2} proved that
if $f$ is expressed by a 
depth-$c$ circuit of size $N$ then 
\begin {equation}
\label {eq:bop}
I(f) \le C_1 \log ^{c-1} N.
\end {equation}

Earlier, Linial, Mansour and Nisan \cite {LMN}
proved that the Fourier 
coefficients of functions which can be expressed by Boolean circuits
of polynomial (or quasi-polynomial) size and bounded depth in AC0 decays 
exponentially above poly-logarithmic ``frequencies''. 
Both these results rely on the fundamental 
H\aa{}stad Switching Lemma, see \cite {Ha,AS}. 

Recall that a {\bf monotone circuit} is one where all  the gates
are monotone increasing in the inputs; i.e., there are no ``not''
gates. The H\aa{}stad lemma for monotone boolean circuits is easier
and was proved already by Boppana \cite {Bop}.

We conjecture that a reverse relation to \ref {eq:bop} 
also holds. 

\begin {conj}[Reverse H\aa{}stad]
\label{RH.conj}
For every $\eps>0$ 
there is a $K=K(\eps )>0$ satisfying the following.
For every monotone $\A\subset\Omega_n$, there is a $\Cal B\subset\Omega_n$
such that $\P[\A\triangle\Cal B]<\eps$ and $\Cal B$ can be expressed
as a Boolean circuit such that
$$
(\log N)^{c-1} < K I(\A),
$$
where $c$ and $N$ are the depth and size of the circuit, respectively.

\end {conj}

Monotone Boolean functions with bounded influence were 
characterized by Friedgut \cite {Fr1,Fr2}. The results of 
\cite {BK} are also relevant to this conjecture.

Ha Van Vu raised the question if there is a spectral way to distinguish 
between bounded depth circuits of polynomial size and bounded depth 
circuits of quasi-polynomial size.  In particular, he was looking for a way 
to show that the graph property ``having a clique of size $\log n$"
for graphs with $n$ vertices, cannot be expressed by a bounded depth 
circuit of polynomial size. (Here the set of variables
correspond to the $n \choose 2$ possible edges.)

\begin {conj}
Let $\epsilon >0$ be a fixed real number and 
$c \ge 1 $ be a fixed integer. 
Let $\A$ be a monotone property expressed by a depth-$c$ 
circuit of size $M$ and let $f=\chi_\A$. Then 
there is a set $\cal S$ of polynomial size in $M$ (where
the polynomial depends on $c$ and $\eps$) so that 
$$\sum \{\hat f^2(S): S \notin {\cal S}\} \le \epsilon. $$ 
\end {conj}

This conjecture may also apply to TC0, see below. It would be of 
great interest to characterize Boolean functions for which most of the
weight
of the Fourier coefficients is concentrated on a set of polynomial size
in $n$.

\subsection {TC0 and noise sensitivity}
 
Noise sensitivity seems related to another class of 
boolean functions - threshold circuits of bounded depths
see \cite {Yao, Has2}. In a threshold circuit each gate 
is a weighted majority function. For the study of spectral properties
of signs of low degree polynomials see Bruck \cite {Bruck} and 
Bruck and Smolensky, \cite {Bruck2}.

\begin {conj}
\label{tco.conj}
 Let $f$ be a boolean function given by a {monotone} threshold circuit
of depth $c$ and size $M$. Then 
\begin {equation}
\label {e:tc0}
J(f)=O(1) (\log M)^{c-1}.
\end {equation}
\end{conj}

Thus,  for $1/\epsilon \le O(1)(\log M)^{c-1}$ 
we expect that $\var(f, \epsilon)$ is bounded away from zero.
Also here it is a tempting conjecture that a reverse
relation holds. 

We conjecture further that all functions 
$f$ that can be expressed by a depth-$c$ 
monotone threshold circuit 
where all the threshold gates are balanced are uniformly stable.
(And in particular, $J(f)= O(1)$.)
Possibly, functions in this class of functions approximate 
arbitrary well arbitrary uniform stable monotone Boolean functions.

Conjecture \ref {tco.conj}
implies theorems of Yau \cite {Yao}
and H\aa{}stad and Goldmann \cite {Has2}. They proved that the and/or tree
(or equivalently the example of ternary tree of Section~\ref {s.ex}) 
does not belong to monotone TC0; i.e., it cannot be expressed 
as a monotone bounded depth circuit of 
polynomial size. 

The results of Yau and H\aa{}stad are still open 
for the non-monotone case. This would follow 
if relation \ref {e:tc0} holds even for
every monotone boolean 
function $f$ given by a (general) threshold circuit
of depth $c$ and size $M$.

\section {Random walks}
\label {s.rw}

For nonempty $\A \subset \Omega_n$, consider a random walk defined as
follows: 
start with a point chosen at random uniformly from $\A$, and at each
step, stay where you are with probability $1/2$, and with probability
$1/(2n)$ move to any one of the neighboring vertices.
Let $\P_{\A}^t$ be the measure on $\Omega_n$ given by the location
of the walk after $t$ steps, and set 
$$
W(\A,\eps)  := \inf\bigl\{t : \|\P_{\A}^t-\P\|<\eps\bigr\}
.
$$
Here $\|\P_{\A}^t-\P\|$ is the measure ($L^1$) norm of the 
difference between $\P_{\A}^t$ and the uniform measure.
 
\begin {thm}
Suppose that $\A_m\subset\Omega_{n_m}$ is
a sequence of events satisfying $\inf_m\P[\A_m]>0$.
\begin{enumerate}
\item $\{\A_m\}$ are asymptotically noise sensitive iff\/
$\lim_m W(\A,\eps)/n_m =0$ for every fixed $\eps>0$.
\item If $\beta (\A_m) \to \beta$,
then $W(\A_m, \epsilon ) \leq n^{1-\beta -o(1)}$.
\end{enumerate}
\end {thm}

\proof
Set $f_t(x):=2^{n_m}\P_{\A_m}^t[\{x\}]$.
Note that
$$
f_{t+1} = (1/2)f_t + (2n_m)^{-1}\sum_{j=1}^{n_m}
\sigma_j f_t.
$$
Consequently,
$$
\hat f_{t}(s) = \left(\frac{2n_m-|s|}{2n_m}\right)^t \hat f_0 (s)
=
\P[\A_m]^{-1}\left(\frac{2n_m-|s|}{2n_m}\right)^t \hat \chi_{\A_m} (s)
.
$$
This gives for every $k=1,2,\dots$
\begin{align*}
\|\P_{\A_m}^t-\P\|^2
&\leq \|f_t-1\|_2^2
=
\P[\A_m]^{-1}\sum_{0\neq s\in\Omega_{n_m}}
\left(\frac{2n_m-|s|}{2n_m}\right)^{2t} \hat \chi_{\A_m} (s)^2
\\&
\leq
\P[\A_m]^{-1}\left(\frac{2n_m-k}{2n_m}\right)^{2t}
+
\sum_{0<|s|<k} \hat \chi_{\A_m} (s)^2
\\&
\leq
\P[\A_m]^{-1}\exp(-tk/n_m)
+
\sum_{0<|s|<k} \hat \chi_{\A_m} (s)^2
.
\end{align*}
The theorem follows.
\ep\medskip

\section{Changing a fixed size set of bits} \label{s.fix}

The noise operator $\NN_\eps$ changes every input variable independently
of the others, and the expected number of bits changed is $\eps n$,
where $n$ is the number of variables. Understanding the effect of 
different types of noise may be of interest.
We consider a variant where a fixed number of bits are changed.
In other words, for $x\in\Omega_n$ and $q\in\nn$, let
$\tilde\NN_q(x)=x\oplus s$, where $s$ is chosen randomly uniformly
among $s\in\Omega_n$ with cardinality $q$, independent
from $x$.  Here $\oplus$ is addition mod $2$; that is, xor.

The analysis of the noise $\tilde \NN_q$ is similar to
that of $\NN_\eps$, but a little care is needed.
Consider the following example.  Let $\Cal P\subset\Omega_n$ 
consist of those $x\in\Omega_n$ such that $|x|$ is odd.
This event $\Cal P$ is called {\bf parity}. 
Observe that for each fixed $q$,
the conditioned probability
$ \P[\tilde \NN_q(x) \in\Cal P| x=y ] $ is either zero or $1$.
In other words, knowing $x$ allows a perfect prediction for
$\tilde\NN_q(x)\in\Cal P$.
Note that $\hat{\chi}_{\Cal P}(S)$ is nonzero only when
$S\in\{\emptyset,\nn\}$.
This means that the vanishing of the weight of the lower
Fourier coefficients does not imply sensitivity to $\tilde\NN_q$,
as in Theorem~\ref{t.foursens}.

For $f:\Omega_n\to \R$ and $q\in\nn$ set
$$
\tilde{\var}(f,q)
= \mbox{var}_y\Bigl( \E[f(\tilde \NN_q(x)) | x=y ] \Bigr)
= \E_y \E\Bigl(f (\tilde \NN_q(x))| x=y \Bigr)^2 - (\E f)^2
.
$$
We say that a sequence of events $\A_m\subset\Omega_{n_m}$
is {\bf asymptotically noise sensitive with respect to $\tilde\NN$}
if for every $\eps\in(0,1)$ and every sequence $\{q_m\}$ with
$\eps n_m \leq q_m \leq (1-\eps)n_m$, we have
$$
\lim_m \tilde{\var}(\A_m,q_m)=0
.
$$

Note that this is equivalent to the straightforward
analog of the definition for asymptotic noise sensitivity
to our current setting.

\begin{thm}\label{t.tilde}
Let $\A_m\subset\Omega_{n_m}$ be a sequence of events,
and set $g_m = \chi_{\A_m}$.
\begin{enumerate}
\item\label{i.i1}
This sequence is
asymptotically noise sensitive with respect to $\tilde\NN$
iff for every finite $k$
$$
\lim_m
\sum\Bigl\{\hat g_m(S)^2 : S\subset\nn, 1\leq |S|\leq k\mbox{ or }
|S|\geq n-k  \Bigr\} = 0
.
$$
\item
A sufficient condition for asymptotic noise sensitivity is
$\H(\A_m)\to 0$.
\end{enumerate}
\end{thm}

\proof
For $f:\Omega_n\to\R$ set
$$
\tilde T_qf(y) = \E f(\tilde \NN_q(y))
.
$$
We now compute the Fourier coefficients of $\tilde T_q f$.
Take $r\in\Omega_n$.
\begin{eqnarray*}
\E(\tilde T_q f\cdot u_r)
& = &
2^{-n}\sum_x \tilde T_qf(x)(-1)^{|r\cap x|}
\\& = &
2^{-n}{n\choose q}^{-1}\sum_x\sum_{|s|=q}
f(x\oplus s)(-1)^{|r\cap x|}
\\& = &
2^{-n}{n\choose q}^{-1}\sum_y\sum_{|s|=q}
f(y)(-1)^{|r\cap y|}(-1)^{|r\cap s|}
\\& = &
\hat f(r)
{n\choose q}^{-1}
\sum_j (-1)^j{|r|\choose j}{n-|r|\choose q-j}
.
\end{eqnarray*}
Consequently,
\begin{equation}\label{e.ttf}
\tilde T_q f = \sum_{r\in\Omega_n} c(n,q,|r|)\hat f(r) u_r,
\end{equation}
where
$$
c(n,q,k) =
{n\choose q}^{-1}
\sum_j (-1)^j{k\choose j}{n-k\choose q-j}
.
$$
Since $c(n,q,0)=1$, this gives,
\begin{equation}\label{e.ttff}
\tilde{\var}(f,q)
= \|T_q f\|_2^2 - \hat f(\emptyset)^2
=\sum_{\emptyset\neq S\subset\nn}
c(n,q,|S|)^2 \hat f(S)^2
.
\end{equation}
Consequently, for~\ref{t.tilde}.\ref{i.i1} it is enough to understand the
behavior of the coefficients $c(n,q,k)$.
For this, consider the sequence  
$$
a_j= {k\choose j}{n-k\choose q-j}{n\choose q}^{-1}
.
$$
The sequence has a unique maximum, which occurs
when $j$ is an integer $j'$ close to $qk/n$.
Consequently, $c(n,q,k)\leq 2 a_{j'}$.
Now let $n,k,q\to \infty$, and assume that
$\eps n\leq q\leq (1-\eps)n$ and $n-k\to\infty$,
where $\eps>0$ is fixed.
Then
$$
\lim  c(n,q,k) = 0.
$$
This gives one direction in~\ref{t.tilde}.\ref{i.i1}.

Also note that when $q<n/(3k)$,
$|c(n,q,k)|$ is approximately $a_0$.
This gives $\liminf |c(n,q,k)|>0$
when $k$ is fixed, $n\to\infty$ and
$q$ is about $n/3k$.
Since $c(n,q,k)=\pm c(n,q,n-k)$,
we get the other direction of~\ref{i.i1}.

Now assume that $\H(\A_m)\to 0$.
{}From Theorem~\ref{t.foursens} we know that
$$
\lim_m
\sum\Bigl\{\hat g_m(S)^2 : S\subset\nn, 1\leq |S|\leq k \Bigr\} = 0
$$
for every fixed $k$.
Equation~(\ref{e.ii}) gives,
\begin{eqnarray*}
\H(\A_m)
&\geq&
 {n_m\over 2}\sum \Bigl\{\hat g_m(S)^2 \hat g_m(S')^2 :
S,S'\subset\nn,\ |S|,|S'|\geq 3 n_m/4\Bigr\}
\\&\geq&
\left(\sum_{|S|\geq 3 n_m/4} \hat g_m(S)^2\right)^2
.
\end{eqnarray*}
Consequently,
$$
\lim_m
\sum\Bigl\{\hat g_m(S)^2 : S\subset\nn, |S|\geq n- k \Bigr\} = 0
$$
and the proof is complete.
\ep

\bigskip
\filbreak
\begingroup
{
\small
\begin{sc}
\parindent=0pt\baselineskip=12pt
\def\email#1{\par\qquad {\tt #1} \smallskip}
\def\emailwww#1#2{\par\qquad {\tt #1}\par\qquad {\tt #2}\medskip}

The Weizmann Institute of Science,
Rehovot 76100, Israel
\emailwww{itai@wisdom.weizmann.ac.il}
{http://www.wisdom.weizmann.ac.il/$\sim$itai/}

The Hebrew University, 
Givat Ram, Jerusalem 91904,
Israel
\emailwww{kalai@math.huji.ac.il}
{http://www.ma.huji.ac.il/$\sim$kalai/}

The Weizmann Institute of Science,
Rehovot 76100, Israel
\emailwww{schramm@wisdom.weizmann.ac.il}
{http://www.wisdom.weizmann.ac.il/$\sim$schramm/}
\end{sc}
}
\filbreak

\endgroup


\begin{thebibliography}{MM}

\let\bbii=\bibitem
\def\bibitem[#1]#2{\bbii{#2}}

\bibitem [ADJ]{ADJ} J. Ambjorn, B. Durhuus and T. Jonsson, {\it Quantum
Geometry},
Cambridge University Press, Cambridge, 1997.

\bibitem [AS]{AS} N. Alon and J. Spencer, {\it The Probabilistic Method}, 
Wiley, New York (1992).

\bibitem [Be]{Be} W. Beckner, Inequalities in Fourier analysis, {\it 
Annals of Math.} {\bf 102} (1975), 159--182.

\bibitem [BL]{BL} M. Ben-Or and N. Linial, Collective 
coin flipping, in {\it Randomness and Computation} (S. Micali, ed.),
Academic Press, New York, (1990), pp. 91--115. Earlier version:
Collective coin flipping, robust voting games, 
and minima of Banzhaf value, Proc. 26th IEEE Symp. on the 
Foundation of Computer Science, (1985), pp. 408--416.

\bibitem[BS]{BS}
I. Benjamini and O. Schramm, 
Conformal invariance of Voronoi percolation, {\it Commun. Math. Phys.}, 
{\bf 197} (1998), 75--107.

\bibitem [BKS2]{BKS2}
I. Benjamini, G. Kalai and O. Schramm, 
Noise sensitivity, concentration of measure and first passage percolation, 
in preparation.
 
\bibitem [Bo]{Bo} A. Bonami, Etude des coefficients Fourier des fonctiones 
de $L^p(G)$, {\it Ann. Inst. Fourier} {\bf 20} (1970), 335--402. 

\bibitem [Bop]{Bop} R. Boppana, Threshold functions and bounded depth
monotone circuits, {\it Proceedings of 16th Annual ACM
Symposium on Theory of Computing} (1984), 475--479. 

\bibitem [Bop2]{Bop2} R. Boppana, The average sensitivity of bounded 
depth circuits, {\it Inform. Process. Lett.} {\bf 63} (1997) 257--261. 


\bibitem[BKKKL]{BKKKL} J. Bourgain, J. Kahn, G. Kalai, Y. Katznelson 
and N. Linial,
The influence of variables in product spaces, {\it Isr. J. Math.} {\bf 77} 
(1992), 55--64.

\bibitem[BK]{BK} J. Bourgain and G. Kalai, Influences of variables and 
threshold intervals under group symmetries, 
{\it Geom. Funct. Anal.}, {\bf 7} (1997), 438-461.

\bibitem [Br]{Bruck} J. Bruck, 
Harmonic analysis of polynomial threshold 
functions. SIAM J. Discrete Math. {\bf 3} (1990), 168--177.  

\bibitem [BrSm] {Bruck2} J. Bruck and R. Smolensky,  
Polynomial threshold functions, ${\rm AC}\sp 0$ functions, and spectral
norms. SIAM J. Comput. {\bf 21} (1992), 33--42. 


\bibitem [BH] {BH} A. Bunde and S. Havlin (ed.s'), {\it Fractals and 
Disordered Systems}, Springer 1991.


\bibitem[CCFS]{CCFS} J. T. Chayes, L. Chayes D. S. Fisher and T. Spencer,
Finite-size scaling and correlation length for disordered systems,
{\it Phys.\ Rev.\ Lett.{}} {\bf 57} (1986), 2999--3002.


\bibitem [Fr1]{Fr1} E.\ Friedgut, Boolean functions with low average 
sensitivity, {\it Combinatorica} {\bf 18} (1998), 27--36.

\bibitem [Fr2]{Fr2} E.\ Friedgut, 
Necessary and sufficient conditions for 
sharp thresholds of graphs properties and the $k$-sat problem, 
{\it Jour. Amer. Math. Soc.} {\bf 12} (1999), 1017--1054.

\bibitem[FK]{FK} E.\ Friedgut and G.\ Kalai,
Every monotone graph property has a sharp threshold, 
{\it Proc.\ Amer.\ Math.\ Soc.} {\bf 124} (1996), 2993--3002.

\bibitem[G]{G} G. Grimmett, {\it Percolation}, Springer-Verlag, Berlin
(1989).

\bibitem[hag96]{hag96} O. H\"{a}ggstr\"{o}m, Y. Peres, and J. E. Steif,
Dynamical percolation, {\it Ann. IHP} {\bf 33} (1997), 497--528.



\bibitem [Ha] {Ha}J. H\aa{}stad, 
Almost optimal lower bounds for small depth circuits,
in {\it Randomness and Computation}, 
{\bf 5}, ed. S. Micali, (1989), 143--170.

\bibitem [Has2] {Has2} J. H\aa{}stad 
and M. Goldmann, On the power of small-depth threshold
circuits, {\it Computational Complexity}, {\bf 1} (1991), 113--129.



\bibitem
[KKL]{KKL} J. Kahn, G. Kalai and N. Linial, The influence of variables
on boolean functions, {\it Proc. 29-th Ann. Symp. on Foundations of 
Comp. Sci.}, (1988), 68--80. 

\bibitem [Ke]{Ke} H. Kesten, Scaling relations for $2D$-percolation,
{\it Comm. Math. Phys.} {\bf 109} (1987), 109--156.

\bibitem [KZ] {KZ} H. Kesten  and Y. Zhang, 
Strict inequalites for some critical exponents in 2D-percolation.
{\em J. Statist.\ Phys.} (1987) {\bf 46}, 1031--1055.

\bibitem [LPS]{LPS} R. P. Langlands, P. Pouliot, and Y. Saint-Aubin 
Conformal invariance in two-dimensional percolation,
{\it Bull. Amer. Math. Soc. (N.S.)} {\bf 30} (1994), 1--61.

\bibitem [LMN]{LMN} N. Linial, Y. Mansour and N. Nisan, Constant depth 
circuits, Fourier transform, and learnability, 
{\it J. Assoc. Comput. Mach.} {\bf 40} (1993), 607--620. 


\bibitem[Petrov]{Petrov} V.~V.~Petrov,
{\it Limit theorems of probability theory},
{Oxford University Press}, {(1995)}.
    

\bibitem [R]{R} L. Russo, A note on percolation, {\it ZW.} {\bf 43} (1978),
39--48.

\bibitem [SW]{SW} P. Seymour and D. Welsh,  Percolation probabilities on the
square lattice.  Advances in Graph Theory.
{\it Ann. Discrete Math.} {\bf 3} (1978), 227--245.

\bibitem [T1]{T1} M. Talagrand, On Russo's approximate zero-one law, 
{\it Ann. of Prob.} {\bf 22} (1994), 1576--1587.

\bibitem [T1.5]{T1.5} M. Talagrand, Concentration of measure and
isoperimetric 
inequalities in product spaces, {\it  Publ. I.H.E.S.}, {\bf 81} (1995),
73--205.



\bibitem [T2]{T2} M. Talagrand,  How much are increasing sets  
positively correlated? {\it Combinatorica} {\bf 16} (1996), 243--258.


\bibitem [Ts1]{Ts1} B. Tsirelson, 
Fourier-Walsh coefficients for a coalescing flow (discrete time), preprint,
math.PR/9903068.  


\bibitem [Ts2]{Ts2} B. Tsirelson, The Five noises, preprint.

\bibitem [Yao]{Yao} A. Yao, Circuits and local computation, 
{\it Proceedings of 21st Annual ACM Symposium on Theory of Computing},
(1989), 
186--196.



\end{thebibliography}
\end{document}